\documentclass[a4paper,10pt,oneside]{amsart}
\pdfoutput=1
\usepackage[utf8]{inputenc}
\usepackage[english]{babel}
\usepackage[T1]{fontenc}
\usepackage{amsmath,amssymb}
\usepackage{subfig}
\usepackage{geometry}
\usepackage{bm}
\usepackage{color}
\usepackage{graphicx}
\usepackage{float}
\usepackage[foot]{amsaddr}
\usepackage{booktabs}
\usepackage{calligra}
\usepackage{comment}
\usepackage{algpseudocode,algorithm,algorithmicx}
\usepackage{bmpsize}

\newcommand{\R}{\mathbb{R}}
\newcommand{\N}{\mathbb{N}}
\newcommand{\D}{\mathbb{D}}

\DeclareMathOperator*{\argmax}{arg\,max}

\graphicspath{{Figures/}}

\newtheorem{remark}{Remark}[section]

\begin{document}

\title[HiMod techniques for flow modeling in a parametrized setting]{Hierarchical model reduction techniques for flow modeling in a parametrized setting}
\author{Matteo Zancanaro$^1$}
\author{Francesco Ballarin$^1$}
\author{Simona Perotto$^2$}
\author{Gianluigi Rozza$^1$}
\address{$^1$ mathLab, Mathematics Area, SISSA, Via Bonomea 265, I-34136 Trieste, Italy, $^2$ MOX, Dipartimento di Matematica, Politecnico di Milano, Piazza L. da Vinci 32, I-20133 Milano, Italy}
\date{}

\begin{abstract}
In this work we focus on two different methods to deal with parametrized partial differential equations in an efficient and accurate way.
Starting from high fidelity approximations built via the hierarchical model reduction discretization, we consider two approaches, both based on a
projection model reduction technique. The two methods differ for the algorithm employed during the construction of the reduced basis.
In particular, the former employs the proper orthogonal decomposition, while the latter relies on a greedy algorithm according to the certified reduced basis technique. 
The two approaches are preliminarily compared on two-dimensional scalar and vector test cases.
\end{abstract}

\maketitle

\section{Introduction}
The interest in fluid dynamics simulations is growing more and more in the scientific community and in society, in terms of both spread and relevance. This is, most of all, due to practical issues and time reasons.
Physical experiments are often very expensive and time demanding so that, for some specific applications (e.g., in naval or aeronautic applications as well as 
in medical surgery planning), they are not well-suited, and numerical simulations become the actual tool for modeling reliable scenarios in such contexts.

Although computational power is continuously growing, 
standard methods in computational fluid dynamics (CFD), such as 
direct numerical simulations based on finite elements, may be very
demanding in terms of computational time and numerical sources, especially
when one is interested in simulating challenging phenomena in complex domains with a certain accuracy, or, 
even more, when dealing with multiquery or parametric frameworks \cite{Lassila2014}.

For these reasons, many different methods have been proposed in the scientific panorama with the aim of offering a compromise between modeling accuracy and computational efficiency. Model reduction techniques represent a relevant solution in such a direction \cite{Lassila2014}.
Some of them are strictly intertwined with the model of the phenomenon at hand, while others perform the reduction only under specific physical assumptions on the described configuration \cite{quarteroni2014reduced,benner2017model,BennerSIAM}.

In this work, we focus on the hierarchical model (HiMod) reduction technique~\cite{ern2008hierarchical,perotto2010hierarchical,perotto2014survey,PerottoZilio13}. 
This procedure has been devised to describe CFD configurations where a principal dynamics overwhelms the transverse ones, with a strong interest for hemodynamic 
configurations~\cite{BarbosaPerotto19,Guzzetti18,Rusconi17}.
The leading dynamics is aligned with the main stream of the flow, while transverse dynamics are generally induced by geometric irregularities in the computational domain and play a role only in localized areas.
In practice, the idea is to discretize the different dynamics by resorting to different 
numerical methods, in the spirit of a separation of variable. 
For instance, in the original proposal of HiMod reduction, the main direction of the flux is discretized with one-dimensional (1D) finite elements, while the transverse dynamics are reconstructed by using few degrees of freedom via a suitable modal basis.
This separate discretization, independently of the dimension of the (full) problem at hand, leads to solving
a system of coupled 1D problems, whose coefficients include the effect of the transverse dynamics.
This ensures the HiMod reduction has a reliability which is considerably higher compared with standard 1D reduced models, and at a computational cost which remains absolutely affordable. 
The computational advantages provided by a HiMod discretization have also been verified for simulations in real geometries~\cite{Aletti18,Guzzetti18,Pablo}. In particular, 
HiMod reduction guarantees a linear dependence of the computational cost on the number of degrees of freedom in contrast to a standard full finite element (FE) model which demands a suitable power of such a number.

The interest of this paper is a parametric setting, where the reference model, coinciding with a
parametric partial differential equation, has to be solved several times, for many different values of the parameter.
The goal we pursue is to approximate, for each value of the parameter, the HiMod discretization
by a modeling procedure which turns out to be computationally cheaper than HiMod reduction itself. \\
A first effort in such a direction is proposed in~\cite{Baroli,LupoPasiniPerottoVeneziani}. The authors
apply a proper orthogonal decomposition (POD) procedure to HiMod approximations,   
to extract a reduced basis which allows us to predict the HiMod discretization associated with any value 
of the parameter. The new procedure, named HiPOD, is numerically investigated on scalar elliptic problems
and on the Stokes equations in~\cite{Baroli}.
In this paper, we investigate a new procedure, alternative to HiPOD, to pursue the same goal
of managing, in a cheap way, a parametric framework. 
In particular, we aim at exploiting the computational advantages provided by a greedy algorithm
in the construction of a reduced basis~\cite{devore2013greedy,hesthaven2014efficient}. 
For this purpose, we combine HiMod reduction with the reduced basis (RB) approach~\cite{rozza2014fundamentals, Rozzabook}, into the new technique called HiRB.
HiPOD or HiRB approximations considerably decrease the computational effort due to the lower dimension of the high fidelity problems. According to an offline/online paradigm, the offline stage
remains the bottleneck from a practical viewpoint. However, the employment of HiMod discretizations as high fidelity 
solutions significantly reduces the computational effort of this phase. 
Finally, a system of very small order is solved during the online phase and yields a reliable  approximation for the parametric problem at hand. 

The paper is organized as follows. Section~\ref{HiModsec} introduces the HiMod setting and particularizes such a procedure both to a scalar advection-diffusion-reaction (ADR) problem and to the Stokes equations. 
Sections~\ref{HiPOD_sec} and~\ref{HiRBSec} exemplify the HiPOD and the HiRB procedures, respectively, on the test problems introduced in the previous section. Particular care is devoted to the inf-sup condition characterizing the discretization
of the Stokes problem. Actually, 
although the high fidelity solutions satisfy the Ladyzhenskaya-Brezzi-Babu{\v s}ka (LBB) condition, this
is not ensured either by the POD or the RB formulation. To overcome this issue, we propose here to 
resort to the supremizer enrichment stabilization technique~\cite{ballarin2015supremizer}.
Section~\ref{compariamo} performs a preliminary comparison between HiPOD and HiRB, starting
from the (2D) test cases considered throughout the paper. Finally, some conclusions are 
drawn in Section~\ref{conclusione} and future developments are summarily itemized.

\section{The HiMod setting}\label{HiModsec}
\begin{figure}[tbh]
\centering
\includegraphics[height=0.24\textwidth]{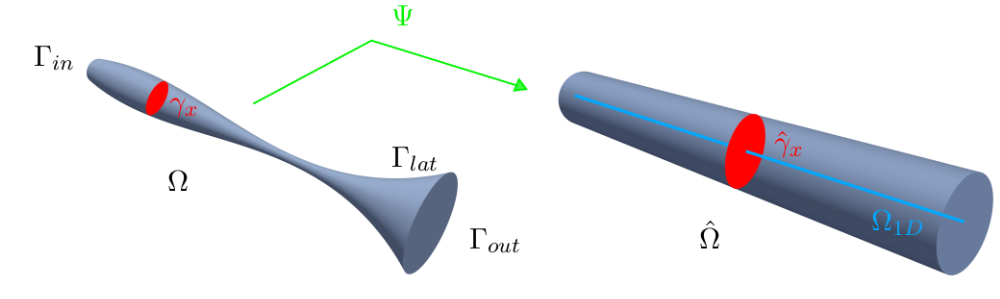}
\caption{HiMod map between the physical and the reference domain.}\label{Map}
\end{figure}

We summarize here the main features of the HiMod reduction technique, following the original setting in 
\cite{ern2008hierarchical,perotto2010hierarchical,PerottoZilio13,perotto2014survey}. 
We assume that the $d$-dimensional domain, 
$\Omega$,  with $d=2, 3$, coincides with the fiber bundle
\begin{equation}
\Omega= \Omega_{1D} \times \gamma_{x},
\label{eq:omega}
\end{equation}
where $\Omega_{1D}$ is the supporting fiber aligned with the main flow, while $\gamma_{x}$ denotes the $(d-1)$-dimensional transverse fiber at point $x\in \Omega_{1D}$, parallel to the secondary dynamics.
In practice, computations are performed in a reference domain, $\hat \Omega$, so that 
$\bm{\Psi}(\Omega)=\hat \Omega$, $\bm{\Psi}: \Omega \to \hat{\Omega}$ being a sufficiently regular map (see Figure~\ref{Map}). 
In general, domain $\hat \Omega$ coincides with a rectangle ($d=2$) or with a right circular cylinder ($d=3$). 
For simplicity, we consider a rectilinear axis $\Omega_{1D}=(0, L)$ with $L > 0$, so that $\bm{\Psi}(x,\bm{y}) = (x, \bm{\psi}_x(\bm{y}))$, for any $(x,\bm{y})\in \Omega$. Map $\bm{\Psi}$ preserves the supporting fiber and deforms only the transverse shape of the domain
via the map $\bm{\psi}_x: \gamma_x \to \hat \gamma$ between the generic, $\gamma_x$ , and the  reference, $\hat \gamma$,
transverse section. This induces a decomposition similar to \eqref{eq:omega} on the reference domain as well, where
$\hat\Omega= \Omega_{1D} \times \hat \gamma$.
We refer to \cite{dd,Rusconi17,BarbosaPerotto19} for the more general case of a curvilinear fiber $\Omega_{1D}$. 

In the next sections, we apply the HiMod discretization to a scalar and to a vector problem, in order to detail the involved procedures.

\subsection{HiMod reduction for advection-diffusion-reaction problems}\label{sec:himodADR}
A generic scalar ADR problem is the \emph{full problem} we are interested in reducing, namely find $u: \Omega \to \mathbb{R}$ such that
\begin{equation}\label{ADRfull}
\begin{cases}
-\nabla \cdot \big(\nu \nabla u\big)(x, \bm{y}) + \bm{b}(x, \bm{y}) \cdot \nabla u(x, \bm{y}) + \sigma(x, \bm{y}) u(x, \bm{y}) = f(x, \bm{y}) & \quad \text{in} \ \Omega,\\[2mm]
u (x, \bm{y}) = g (x, \bm{y})& \quad \text{on} \ \Gamma_{D},\\[2mm]
\nu(x, \bm{y}) \displaystyle\frac{\partial u}{\partial \bm{n}} (x, \bm{y}) = h(x, \bm{y}) & \quad \text{on} \  \Gamma_{N},\\[3mm]
\nu(x, \bm{y}) \displaystyle\frac{\partial u}{\partial \bm{n}}(x, \bm{y}) + \rho(x, \bm{y}) u(x, \bm{y}) = l(x, \bm{y}) & \quad \text{on} \; \Gamma_{R},
\end{cases}
\end{equation}
with $\Gamma_{D}$, $\Gamma_{N}$, $\Gamma_{R}$ portions of the boundary $\partial \Omega$ 
of $\Omega$, such that
$\stackrel{\hspace*{-.3cm}\circ}{\Gamma_D}\cap \stackrel{\hspace*{-.3cm}\circ}{\Gamma_N}\cap \stackrel{\hspace*{-.3cm}\circ}{\Gamma_R}=\emptyset$ and $\Gamma_D\cup \Gamma_N \cup \Gamma_R=\partial \Omega$, $\bm{n}$ being the unit outward normal vector to $\partial \Omega$. Concerning the problem data, $\nu\in L^\infty(\Omega)$, with $\nu(x, \bm{y}) \ge \nu_0>0$ a.e. in $\Omega$, denotes the diffusion coefficient, $\bm{b}=[b_{x} , b_{\bm{y}}]^T\in [ L^\infty(\Omega)]^d$, with $\nabla \cdot \bm{b}\in L^2(\Omega)$, the advective field, $\sigma\in L^2(\Omega)$, with $\sigma(x, \bm{y})\ge 0$ a.e. in $\Omega$, the reaction, $f\in L^2(\Omega)$ the forcing term, $g\in H^{1/2}(\Gamma_D)$, $h\in L^2(\Gamma_N)$ and $l\in L^2(\Gamma_R)$ are the boundary data, with $\rho\in L^\infty(\Gamma_R)$ and
where standard notation are adopted for function spaces~\cite{Ciarlet}.

HiMod reduction applies to the weak form of the full problem,
\begin{equation}\label{FullProblem}
\mbox{find\ } u \in V\ :\  a(u,v) = F(v) \quad \forall v \in V,
\end{equation}
with $V=H^1_{\Gamma_D}(\Omega)$,
\begin{equation}\label{forms}
\begin{array}{lll}
a(u, v) &= &\displaystyle \int_{\Omega} \nu \nabla u \cdot \nabla v \ d\Omega \ + \int_{\Omega} v \bm{b} \cdot \nabla u \ d\Omega \ + \int_{\Omega} \sigma u v \ d\Omega \ + \int_{\Gamma_R} \rho u v \ dS , \\[4mm]
F(v) &=& \displaystyle \int_{\Omega} f v \ d\Omega \ + \int_{\Gamma_R} l v \ dS \ +\int_{\Gamma_N} h v \ dS ,
\end{array}
\end{equation}

where, to simplify notation, we assume $g=0$ in \eqref{ADRfull} and we drop the dependence on $(x, \bm{y})$.
The assumptions above on the problem data ensure the well-posedness of \eqref{FullProblem}~\cite{Ciarlet}.

Thus, the HiMod formulation for problem \eqref{ADRfull} can be stated as 
\begin{equation}\label{HiModProblem}
\mbox{find\ } u_{m} \in V_{m}\ :\  a(u_{m},v_{m}) = F(v_{m}) \quad \forall v_{m} \in V_{m},
\end{equation}
for a certain $m\in \N^+$, and with 
\begin{equation}\label{HiModSpace}
V_{m} = \bigg\{ v_{m} (x,\bm{y}) = \sum_{k=1}^{m} \tilde{v}_{k} (x) \varphi_{k}(\psi_x(\bm{y}))\ \text{with} \ \tilde{v}_{k}
 \in V_{1D}^{h},\  x \in \Omega_{1D},\  \bm{y} \in \gamma_{x} \bigg\},
\end{equation}
the HiMod space, where $V_{1D}^{h}$ is a 1D
discrete subspace of $H^{1}(\Omega_{1D})$ associated with a subdivision, ${\mathcal T}_h$, of $\Omega_{1D}$, 
$\{ \varphi_{k} \}_{k=1}^m$ is a modal basis of functions defined on $\hat \gamma$, 
orthonormal with respect to the $L^2(\hat \gamma)$- scalar product, and 
$m$ is the modal index, i.e., the number of modes employed to model the transverse dynamics. In what follows, we identify $V_{1D}^{h}$ with the standard space of the (continuous) finite elements~\cite{Ciarlet}, 
while referring to~\cite{dd,Rusconi17,BarbosaPerotto19} for different discretizations.
As far as the choice of $m$ is concerned, it can be fixed a priori, thanks to heuristic considerations
or to a partial knowledge of the full problem, or a posteriori, driven by a modeling error analysis as in~\cite{perotto2014coupled,PerottoZilio15}.  

The HiMod space has to be endowed with
a conformity and a spectral approximability assumption to ensure the well-posedness of formulation \eqref{HiModProblem}, and a standard density hypothesis has to be advanced on the discrete space $V_{1D}^{h}$ to guarantee the convergence of the HiMod approximation, $u_{m}$, to $u$ (we refer to~\cite{perotto2010hierarchical} for the details). 

Concerning the boundary conditions completing problem
\eqref{ADRfull}, we have to distinguish between data assigned on the inflow/outflow boundaries and on the lateral surface of $\Omega$.  In the first case, we employ a modal expansion of the data to be imposed in an essential way.
With reference to lateral boundary conditions, we resort to the approach proposed in~\cite{Aletti18}, where
the authors set a general way to incorporate, essentially, the lateral boundary data by defining a customized basis referred to as an
\emph{educated} modal basis. The effectiveness of such a procedure is successfully investigated both from a theoretical and a numerical point of view in the same work. In the numerical assessment below, we resort to an educated modal basis to manage the lateral boundary conditions.

From a computational viewpoint, discretization \eqref{HiModProblem} turns the full model \eqref{ADRfull} into a system of $m$ coupled 1D problems defined on $\Omega_{1D}$. This represents the strength-point of a HiMod formulation due to the expected saving in terms of computational effort, for $m$ reasonably small. Actually, we are led to solve the HiMod linear system 
\begin{equation}\label{HiModProb}
A_{m} {\mathbf{u}}_{m} = \mathbf{f}_{m},
\end{equation}
where $A_{m} \in \R^{mN_{h} \times mN_{h}}$ and ${\mathbf{f}}_{m} \in \R^{mN_{h}}$ are the HiMod stiffness matrix and right-hand side associated with the bilinear and linear forms in \eqref{forms}, with $N_h={\rm dim}(V_{1D}^{h})$, and where
$\mathbf{u}_{m} \in \R^{mN_{h}}$ collects the (unknown) coefficients of the HiMod expansion
\begin{equation}\label{HiModVar}
u_{m} (x,\bm{y}) = \sum_{k=1}^{m} \sum_{i=1}^{N_{h}} \tilde{u}_{k,i} \theta_{i} (x)  \varphi_{k} (\psi_x(\bm{y})) 
\end{equation}
with $\{\theta_{i}\}_{i=1}^{N_{h}}$ the FE basis.
For a full characterization of system \eqref{HiModProb}, we refer to~\cite{ern2008hierarchical,perotto2010hierarchical}.
\begin{figure}[t]
\centering
\includegraphics[height=0.24\textwidth]{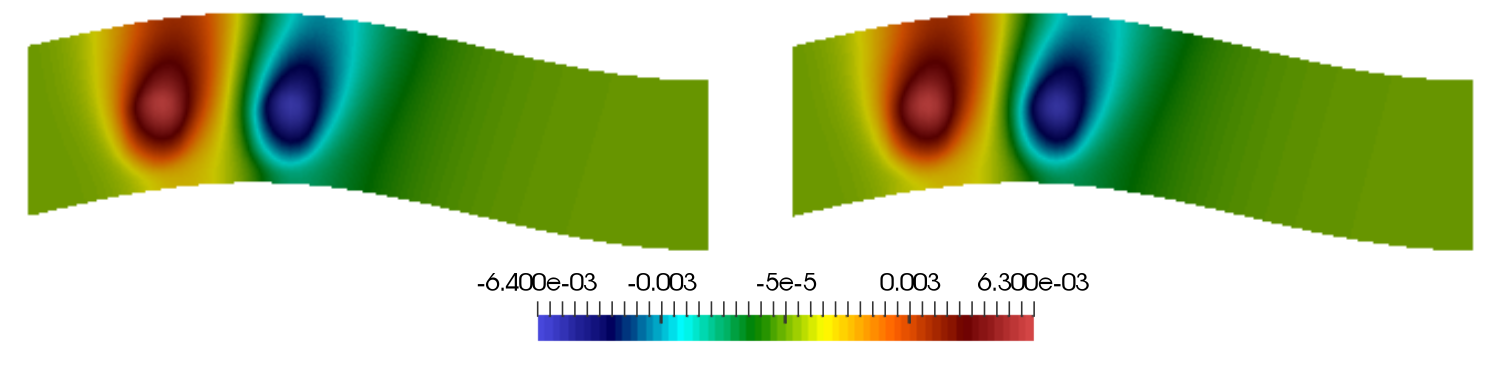}
\caption{ADR test case: comparison between the full solution (left) and the HiMod approximation, $u_8$, (right).}\label{adr_f1}
\end{figure}

To qualitatively investigate the reliability of the HiMod reduction, we solve problem \eqref{ADRfull} 
by means of a FE solver and of a HiMod discretization
on the 2D domain, $\Omega$, identified by the map $\psi_x(y) = y - 0.2 \sin\big(3\pi x/(2L)\big)$, with $x\in [0, 4]$,
$L=4$, and $\hat\Omega = (0, 4) \times (-0.5, 0.5)$.
Concerning the problem data, we assign $\nu=5$, $\bm{b}=[20, 75]^T$, $\sigma=25$, $f=f(x, y) = 1.8\, \chi_{S_1}(x, y)- 1.8\, \chi_{S_2}(x, y)$, with $\chi_\omega$ the characteristic function associated with the subset $\omega\subset \Omega$, 
$S_1=\{(x, y)\ :\ 0.5 (x-0.75)^2 + 0.4 y^2 - 0.02 < 0\}$, $S_2=\{ (x, y)\ :\ 0.5 (x-1.5)^2 + 0.4 y^2 - 0.02 < 0\}$; $\Gamma_D$ and $\Gamma_N$ are identified with the inflow and the outflow boundary, respectively, while $\Gamma_R$ coincides with the lateral surface, where $g=h=l=0$ and $\rho=1$. \\
The FE solver employs affine finite elements on a 2D unstructured mesh consisting of $12800$ triangles.
The HiMod reduction discretizes the main stream with linear finite elements associated with a uniform subdivision of the supporting fiber into $80$ intervals, while resorting to $m=8$ educated modal basis functions in the transverse direction.
Figure~\ref{adr_f1} compares the FE with the HiMod approximation and highlights the good qualitative matching between the two solutions. A quantitative investigation of the HiMod procedure is beyond the goal of this paper and can be found, e.g., in~\cite{perotto2010hierarchical,Aletti18,Guzzetti18}, together with a modeling convergence analysis with respect to both the modal expansion and the FE discretization.

\subsection{HiMod reduction for the Stokes equations}\label{StokesPODProblem}
In this section we generalize the HiMod procedure to a vector problem, namely, to the Stokes equations; find $\bm{u}: \Omega \to \mathbb{R}^d$ and $p: \Omega \to \mathbb{R}$ such that
\begin{equation}
\begin{cases}
-\nabla \cdot\big(2 \nu \D (\bm{u}) \big)(x, \bm{y}) + \nabla p(x, \bm{y}) = \bm{f} (x, \bm{y}) &\quad  \text{in} \ \Omega,\\[1mm]
\nabla \cdot\bm{u} (x, \bm{y}) = 0 & \quad \text{in} \ \Omega,\\[1mm]
u_{y}(x, \bm{y})  =  0  & \quad \text{on} \  \Gamma_{\rm in} \cup \Gamma_{\rm out},\\[1mm]
-\displaystyle\frac{\partial u_{x}}{\partial x}(x, \bm{y}) +  p(x, \bm{y}) =  - C_{\rm in}(x, \bm{y}) & \quad \text{on} \  \Gamma_{\rm in},\\[3mm]
-\displaystyle\frac{\partial u_{x}}{\partial x}(x, \bm{y}) +  p(x, \bm{y}) =  C_{\rm out}(x, \bm{y}) & \quad \text{on} \  \Gamma_{\rm out},\\[3mm]
\bm{u} (x, \bm{y}) =  \bm{0}  &\quad  \text{on} \ \Gamma_{\rm w},\\[1mm]
\end{cases}
\label{PODStokesProblem}
\end{equation}
where $\bm{u}=(u_x, u_{\bm{y}})^T$ and $p$ denote the velocity and the pressure of the flow, 
$\nu>0$ is the kinematic viscosity, $\D (\bm{u})=0.5\big(\nabla \bm{u} + (\nabla \bm{u})^T\big)$ is the strain rate, 
$\bm{f}=[f_x, f_{\bm{y}}]^T$ is the force per unit mass, and $C_{\rm in}$ and $C_{\rm out}$ are the inflow and outflow data, respectively.
The boundary $\partial\Omega$ is partitioned so that the inflow and outflow sections, $\Gamma_{\rm in}$ and $\Gamma_{\rm out}$, coincide with the fibers $\gamma_{0}$ and $\gamma_{L}$, respectively, while $\Gamma_{\rm w}$ denotes the lateral walls $\Omega_{1D} \times \partial\gamma_{x}$. On $\Gamma_{\rm in}$ and $\Gamma_{\rm out}$ we impose a nonhomogeneous tangential Neumann condition, with $C_{\rm in}$ and $C_{\rm out}$ constant values, and we assume the transverse component of the velocity to be null. Finally, a no-slip boundary condition is enforced on the velocity along the wall surface.

By introducing the bilinear forms $a(\cdot, \cdot): \bm{V} \times \bm{V} \to \mathbb{R}$ and $b(\cdot, \cdot): \bm{V} \times Q \to \mathbb{R}$, defined as
\begin{equation*}
a(\bm{u}, \bm{v}) = \int_{\Omega} 2 \nu \D (\bm{u}) : \nabla \bm{v} \ d\Omega,  \quad b(\bm{u}, q) = \int_{\Omega} \nabla \cdot \bm{u} \, q \ d\Omega,
\end{equation*}
and the functional $F(\cdot): \bm{V} \to \mathbb{R}$, given by
\begin{equation*}
F(\bm{v}) = \int_{\Omega} \bm{f}  \cdot \bm{v} \ d\Omega \ + \int_{\partial \Omega} C \bm{n} \cdot \bm{v} \ dS,
\end{equation*}
where $C = -C_{\rm in}$ on $\Gamma_{\rm in}$ and $C = C_{\rm out}$ on $\Gamma_{\rm out}$, with $\bm{V} = \{\bm{v} \in H^1(\Omega; \mathbb{R}^d): v_y = 0 \text{ on } \Gamma_{\rm in} \cup \Gamma_{\rm out} \text{ and } \bm{v} = \bm{0} \text{ on } {\Gamma_{\rm w}}\}$ and $Q = L^2(\Omega)$,  the 
weak form of \eqref{PODStokesProblem} can be stated as finding $\bm{u} \in \bm{V}, p \in Q$ such that
\begin{equation}\label{FullProblemStokes1}
\begin{cases}
a(\bm{u}, \bm{v}) + b(\bm{v}, p) = F(\bm{v}) & \forall \bm{v} \in V,\\
b(\bm{u}, q) = 0 & \forall q \in Q,
\end{cases}
\end{equation}
where the notation is simplified by removing the dependence on $(x, \bm{y})$ and where the natural boundary conditions still have to be properly included in $F$.

The generalization of a HiMod reduction to the Stokes equations deserves particular attention, especially with reference to 
the two-field formulation involved by saddle point problems~\cite{brezzi2013mixed,brezzi1974existence}.
While the search for inf-sup (or LBB) compatible spaces for velocity and pressure is largely investigated for standard FE and spectral discretizations~\cite{brezzi2013mixed,brezzi1974existence,CHQZ}, we are not aware of any theoretical result for hybrid methods involving both techniques. In~\cite{Guzzetti18, Pablo}, empirical criteria to select the HiMod velocity and pressure are provided and numerically checked. 
A first theoretical assessment of these criteria is currently under investigation~\cite{Olivia}.\\
The reduced spaces involved in the HiMod discretization of the Stokes equations are
\begin{align*}
\bm{V}_{m_{\bm{u}}} = \left\{\bm{v}_{m_{\bm{u}}}(x,\bm{y}) = (v_{x,m_{\bm{u}}} (x,\bm{y}), v_{{\bm{y}},m_{\bm{u}}} (x,\bm{y}))^T\ : \ v_{x,m_{\bm{u}}}\in V_{m_{\bm{u}}}, v_{{\bm{y}},m_{\bm{u}}} \in [V_{m_{\bm{u}}}]^{d-1} \right\},\\
Q_{m_p} = \Big\{ q_{m_p} (x,\bm{y}) = \sum_{k=1}^{m_p} \tilde{q}_{k} (x) \eta_{k}(\psi_x(\bm{y})),\ \text{with} \ \tilde{q}_{k}
 \in Q_{1D}^{h},\ x \in \Omega_{1D},\ \bm{y} \in \gamma_{x} \Big\},
\end{align*}
for the velocity and the pressure, respectively, where space $V_{m_{\bm{u}}}$ is the scalar space defined as in \eqref{HiModSpace}.
Notice that we employ the same (educated) modal basis, $\{ \varphi_{k} \}_{k=1}^{m_{\bm{u}}}$, for all the components of the HiMod velocity, while we resort to the modal basis $\{ \eta_{k} \}_{k=1}^{m_p}$ to discretize the pressure.
 
Concerning the compatibility of the velocity with the pressure HiMod spaces, in this work we adopt the empirical criterion in~\cite{Aletti18,Guzzetti18}, so that
we set $m_{\bm{u}}=m_p+2$ and we choose the 1D FE pair $(V_{1D}^{h}, Q_{1D}^{h})$ as the Taylor-Hood P2/P1 elements~\cite{brezzi2013mixed}. We denote by $N_{h,\bm{u}}$ and $N_{h,p}$ the dimension of $V_{1D}^{h}$ and $Q_{1D}^{h}$, so that the dimension of $\bm{V}_{m_{\bm{u}}}$ and $Q_{m_p}$ becomes $dm_{\bm{u}} N_{h,\bm{u}}$ and $m_{p} N_{h,p}$, respectively.

Let us now describe the algebraic formulation for the HiMod discretization of the Stokes problem. After assembling the matrices $A_{m_{\bm{u}}}\in \R^{dm_{\bm{u}} N_{h,\bm{ u}}\times dm_{\bm{u}} N_{h,\bm{u}}}$, $B_{m_p,m_{\bf u}}\in \R^{m_{p} N_{h,p}\times dm_{\bm{u}} N_{h,\bm{u}}}$ and the vector $\mathbf{f}_{m_{\bm{u}}}\in \R^{dm_{\bm{u}} N_{h,\bm{u}}}$ associated with the HiMod discretization of the forms $a(\bm{u}, \bm{v})$, $b(\bm{u}, q)$, and $F(\bm{v})$ in \eqref{FullProblemStokes1}, the linear system 
\begin{equation}
\left[\begin{array}{ll} A_{m_{\bm{u}}} & B^T_{m_p,m_{\bf u}}\\ B_{m_p,m_{\bf u}} & 0 \end{array}\right]
\left[\begin{array}{c} \mathbf{u}_{m_{\bm{u}}} \\ \mathbf{p}_{m_p} \end{array}\right] 
= \left[\begin{array}{c} \mathbf{f}_{m_{\bm{u}}} \\ \bm{0} \end{array}\right]
\label{HiModProbStokes}
\end{equation}
has to be solved, where $\mathbf{u}_{m_{\bm{u}}}\in \R^{dm_{\bm{u}} N_{h,\bm{u}}}$ and $\mathbf{p}_{m_p}\in \R^{m_{p} N_{h,p}}$ collect the unknown coefficients of the HiMod expansion for the velocity, $\bm{u}_{m_{\bm{u}}}$, and the pressure, $p_{m_p}$, respectively, and with $\bm{0}$ the null vector in $\R^{m_{p} N_{h,p}}$.

\begin{figure}[tbh]
\centering
\includegraphics[height=0.22\textwidth]{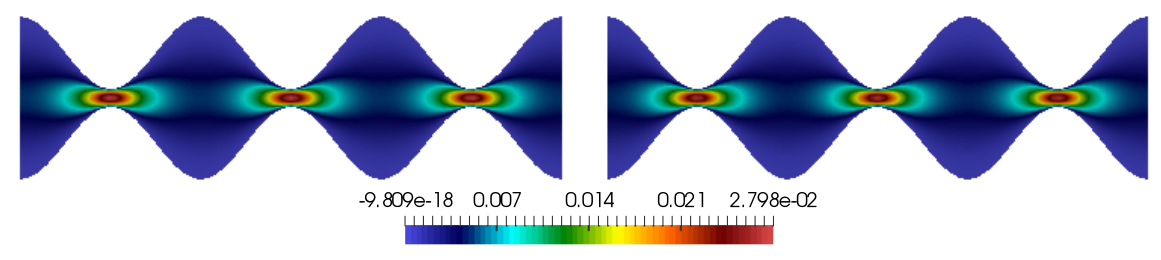}
\caption{Stokes test case: comparison between the full solution (left) and the HiMod approximation, $\bm{u}_7$/$p_5$ (right), for the horizontal component of the velocity.}\label{stokes_f1}
\end{figure}

\begin{figure}[tbh]
\centering
\includegraphics[height=0.24\textwidth]{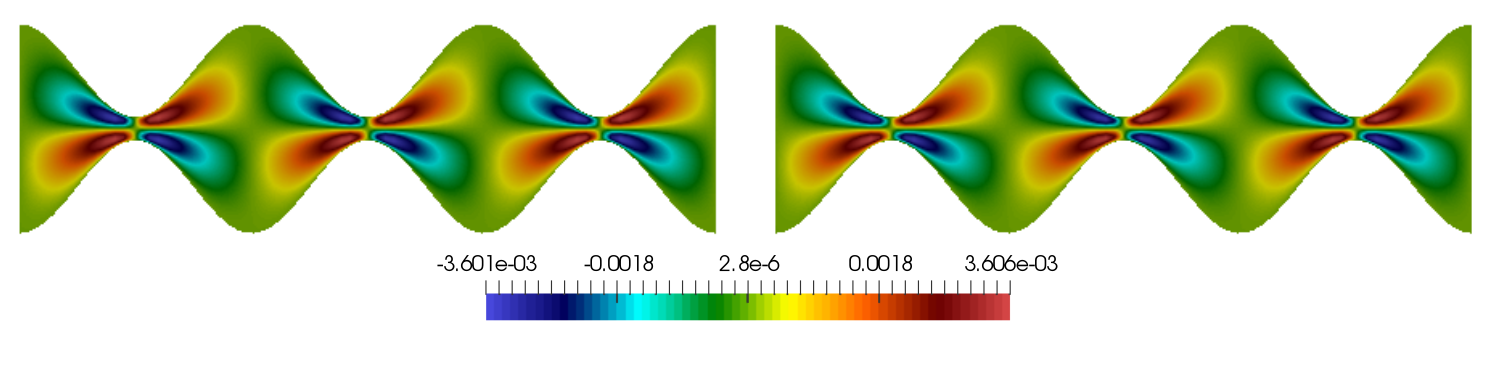}
\caption{Stokes test case: comparison between the full solution (left) and the HiMod approximation, $\bm{u}_7$/$p_5$ (right), for the vertical component of the velocity.}\label{stokes_f2}
\end{figure}

\begin{figure}[tbh]
\centering
\includegraphics[height=0.22\textwidth]{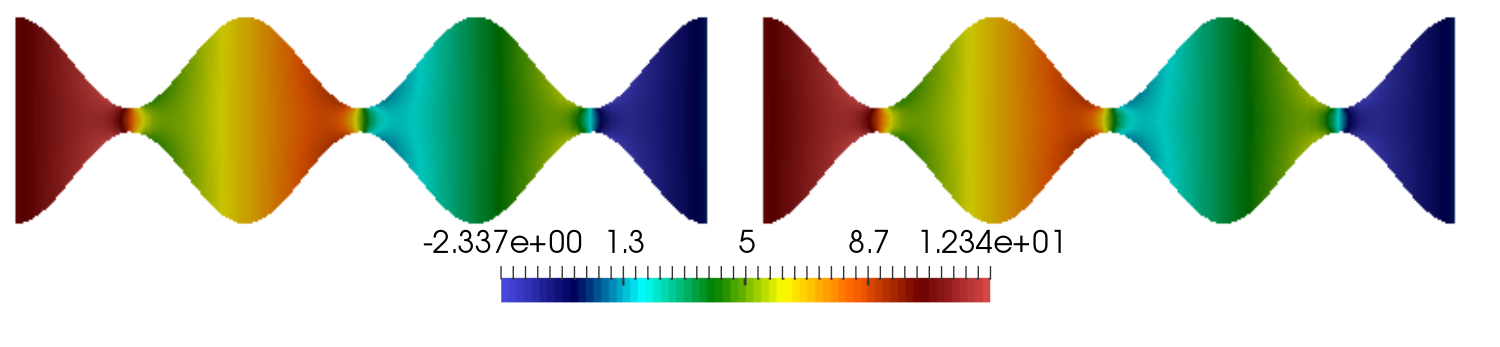}
\caption{Stokes test case: comparison between the full solution (left) and the HiMod approximation, $\bm{u}_7$/$p_5$ (right), for the pressure.}\label{stokes_f3}
\end{figure}

We conclude this section by exemplifying the HiMod procedure on a benchmark Stokes test case. The reference domain is $\hat\Omega = (0, L) \times (-0.5, 0.5) \subset \mathbb{R}^2$, while the map $\psi_x$ is given by
\begin{equation*}
\psi_x(y) = \frac{y}{1 + \frac{2}{5} \sin\Big(\frac{6\pi x}{L}+\frac{\pi}{2}\Big)\frac{2}{H}},
\end{equation*}
so that $\Omega$ coincides with a sinusoidal domain. In particular, we select $L = 6$ and $H = 1$.
Furthermore, we assign $\nu=5$, $\bm{f}=[3, 0]^T$, $C_{in}=10$, $C_{out}=0$.\\
Concerning the HiMod discretization, we enrich the Taylor-Hood P2/P1 discretization of the mainstream by resorting to $m_p = 5$ and $m_{\bm{u}}=7$ 
educated modes to discretize the transverse components of the pressure and the velocity, respectively.
In particular, both the finite element approximations rely on a uniform subdivision of the supporting fiber into $80$ subintervals.
Figures~\ref{stokes_f1}-\ref{stokes_f3} compare the HiMod approximation with a full P2/P1 FE solution computed on an unstructured mesh consisting of $12800$ elements. The two discretizations lead to fully comparable approximations in terms of both velocity and pressure.

\begin{remark}
As investigated in more detail in~\cite{ZancanaroThesis}, the conservative form \eqref{PODStokesProblem}
of the Stokes problem allows one to obtain a more accurate HiMod approximation with respect to a nonconservative formulation. This is due to the coupling between the velocity components (namely, between the off-diagonal blocks of the HiMod matrix) ensured by the conservative form.
\end{remark}

\begin{remark}
Cylindrical domains demand a careful selection of the modal basis as investigated in~\cite{Guzzetti18}, where
a polar coordinate system is employed to model the transverse dynamics in circular and elliptic pipes. As a first alternative, one can resort to the transversally enriched pipe element method (TEPEM)~\cite{Pablo}. In such a case, the physical domain is mapped to a reference slab, so that the modal basis coincides with the
tensor product of the 1D modal functions. In~\cite{Guzzetti18} TEPEM is 
compared with the HiMod reduction based on polar coordinates to highlight pros and cons of the two approaches.
As expected, TEPEM turns out to be easier to implement but less accurate than HiMod, 
especially in the presence of highly oscillatory flows. The higher reliability characterizing the HiMod scheme can be ascribed to the tight
correspondence between the domain geometry and the modal basis. Another alternative to a polar coordinate system is represented by the isogeometric version of the HiMod approach, as recently investigated in~\cite{BarbosaPerotto19} for patient-specific geometries.
\end{remark}

\section{The HiPOD approach}\label{HiPOD_sec}
HiMod reduction allows one to recast a $d$-dimensional problem as a system of 1D problems. Even though this leads to a computational benefit, the overall computational cost still might not be negligible when dealing with multiquery or inverse problems or, more generally, with parametrized settings. In this context, to further lighten the computational effort, we rely on projection-based model reduction techniques, by properly combining the HiMod discretization with the POD. 
In particular, we adopt the standard offline/online paradigm~\cite{Rozzabook}.
The offline phase is meant to build the POD basis, starting from the hierarchical reduction of a certain number of full problems associated with a sampling of the parameter domain. The online phase approximates the HiMod discretization for any new value of the parameter by employing the POD basis. The combination between HiMod and POD justifies the name of this method, HiPOD~\cite{Baroli,LupoPasiniPerottoVeneziani,Meglioli}.

\subsection{HiPOD reduction for ADR problems}
We generalize problem \eqref{FullProblem} to a parameter dependent setting, so that the new problem is
\begin{equation}\label{ParamFullProblem}
\mbox{find\ } u(\bm{\mu}) \in V\ :\   a(u(\bm{\mu}),v; \bm{\mu}) = F(v; \bm{\mu})  \quad \forall v \in V, 
\end{equation}
where $\bm{\mu} \in D\subset \mathbb{R}^P$ denotes a vector of $P$ real numbers collecting the problem parameters, and $D$ is the parameter domain. 

\subsubsection{The offline phase}\label{OfflinePhase}
We choose the sampling $S = \{\bm{\mu}^{(1)}, \bm{\mu}^{(2)}, \ldots, \bm{\mu}^{(M)}\}\subset D^M$ for the parameter $\bm{\mu}$. For each value $\bm{\mu}^{(j)} \in S$, we approximate the corresponding solution, $u(\bm{\mu}_j)$, to \eqref{ParamFullProblem} by computing the HiMod discretization, $u_{m}(\bm{\mu}^{(j)})$,
for a certain value, $m$, of the modal index. According to the modal expansion in \eqref{HiModVar}, this yields the $M$ vectors
\begin{equation*}
{\mathbf{u}}_{m} (\bm{\mu}^{(j)}) = [\tilde{u}_{1,1}^{\bm{\mu}^{(j)}}, \ldots, \tilde{u}_{1,N_h}^{\bm{\mu}^{(j)}}, \ldots, \tilde{u}_{m,1}^{\bm{\mu}^{(j)}}, \ldots, \tilde{u}_{m,N_h}^{\bm{\mu}^{(j)}}]^T \in \R^{mN_{h}},\quad j=1, \ldots, M,
\end{equation*}
collecting, by mode, the HiMod coefficients. These vectors are employed to assemble the response matrix 
\begin{equation*}
U \: = \left[  {\mathbf{u}}_{m} (\bm{\mu}^{(1)}), {\mathbf{u}}_{m} (\bm{\mu}^{(2)}), \ldots, {\mathbf{u}}_{m} (\bm{\mu}^{(M)}) \right]=
\begin{bmatrix}
\tilde{u}_{1, 1}^{\bm{\mu}^{(1)}} & \tilde{u}_{1, 1}^{\bm{\mu}^{(2)}} & \hdots & \tilde{u}_{1,1}^{\bm{\mu}^{(M)}}\\
\vdots & \vdots & \vdots & \vdots\\
\tilde{u}_{1, N_h}^{\bm{\mu}^{(1)}} & \tilde{u}_{1, N_h}^{\bm{\mu}^{(2)}} & \hdots & \tilde{u}_{1,N_{h}}^{\bm{\mu}^{(M)}}\\
\vdots & \vdots & \vdots & \vdots\\
\tilde{u}_{m,1}^{\bm{\mu}^{(1)}} & \tilde{u}_{m,1}^{\bm{\mu}^{(2)}} & \hdots & \tilde{u}_{m,1}^{\bm{\mu}^{(M)}}\\
\vdots & \vdots & \vdots & \vdots\\
\tilde{u}_{m,N_h}^{\bm{\mu}^{(1)}} & \tilde{u}_{m,N_h}^{\bm{\mu}^{(2)}} & \hdots & \tilde{u}_{m,N_{h}}^{\bm{\mu}^{(M)}}\\
\end{bmatrix} \in \R^{mN_{h} \times M},
\end{equation*}
which will be used to extract the POD basis. To this aim, we define 
the correlation matrix associated with $U$, 
\begin{equation}\label{matrixC}
C = U^{T} X_{m,u} \, U \in \R^{M \times M},
\end{equation}
with $X_{m,u} \in \R^{mN_{h} \times mN_h}$ the HiMod matrix associated to the inner product in $V_m$. Then, we consider the spectral decomposition of matrix $C$, so that
\begin{equation}\label{spettrodellemiebrame}
C {\bm \varphi}^*_{k} = \lambda_{k} {\bm \varphi}^*_{k} \quad k = 1, \ldots, M,
\end{equation}
with ${\bm \varphi}^*_{k}/\lambda_{k}$ the $k$th
eigenvector/eigenvalue pair of $C$, where ${\bm \varphi}^*_{k} \in \R^{M}$ and $\lambda_{k}\in \R$. The POD basis is thus identified by the vectors
\begin{equation}
{\bm \varphi}_{k} =  \frac{1}{\lambda_{k}} U {\bm \varphi}^*_{k} \in \R^{mN_{h}}, \quad k = 1, \hdots,  N,
\end{equation}
with $N \le M$. Integer $N$ can be selected driven by heuristic considerations 
(e.g., by studying the trend of the spectrum of $C$) or by an energy criterion; for instance, we pick $N$ such that
\begin{equation*}
E(N)  >  1-\varepsilon \quad \mbox{with}\quad E(N) =  \displaystyle \frac{\sum_{i=1}^{N} \lambda_{i}}{\sum_{i=1}^{M} \lambda_{i}}  
\end{equation*}
and $\varepsilon$ a user-defined tolerance. Independently of the adopted criterion, we denote  
the reduced POD space by $V_{m,N} = \mbox{span} \{ {\bm \varphi}_{1} , \hdots , {\bm \varphi}_{N} \}$, and the matrix collecting, by column, the POD basis functions by $\Phi_{m,N} = [{\bm \varphi}_{1} , \hdots , {\bm \varphi}_{N}] \in \R^{mN_h\times N}$.

\begin{remark}
As an alternative to the approach based on the correlation matrix, one can exploit directly the spectral properties of
the response matrix $U$ to extract the reduced POD basis, by setting $X_{m,u} = I$ in \eqref{matrixC}, with $I \in \R^{mN_{h} \times mN_h}$ the identity matrix~\cite{Baroli}. This twofold possibility is justified by the relation between the singular vectors of $U$ and the eigenvectors of $C$~\cite{Golub}.
\end{remark}

\subsubsection{The online phase}\label{OnlinePhase}
The goal of the online phase is to build a HiMod approximation to problem \eqref{ParamFullProblem}
for any value $\bm{\mu} \in D$ of the parameter, by skipping the solution of the 
associated HiMod system \eqref{HiModProb},
\begin{equation}\label{HiMod_param}
A_{m}(\bm{\mu}) {\mathbf{u}}_{m}(\bm{\mu}) = \mathbf{f}_{m}(\bm{\mu}),
\end{equation}
where the dependence on the parameter $\bm{\mu}$ has been highlighted.
This task is accomplished by means of a projection step, i.e., by solving the system 
\begin{equation}\label{ProjPro}
A_{m, N} (\bm{\mu}) {\mathbf{u}}_{m, N}(\bm{\mu}) = \mathbf{f}_{m, N} (\bm{\mu}),
\end{equation}
with ${\mathbf{u}}_{m, N}(\bm{\mu}) \in \R^N$,
\begin{equation*}
A_{m, N} (\bm{\mu}) = \Phi_{m, N}^{T} A_{m} (\bm{\mu}) \Phi_{m, N} \in \R^{N \times N}, \qquad \mathbf{f}_{m, N} (\bm{\mu}) = \Phi_{m, N}^{T} \mathbf{f}_{m} (\bm{\mu}) \in \R^N.
\end{equation*}
Notice that the order of system \eqref{ProjPro} is significantly smaller compared with 
the HiMod system \eqref{HiMod_param}, where, in general, $N\ll mN_h$. Successively, ${\mathbf{u}}_{m, N}(\bm{\mu})$ is projected back to the original HiMod space, thus obtaining the approximation
\begin{equation*}
{\mathbf{u}}_{m}(\bm{\mu}) \approx \Phi_{m, N} {\mathbf{u}}_{m, N}(\bm{\mu}):= {\mathbf{u}}_{m, N, {\rm POD}}(\bm{\mu}).
\end{equation*} 
In what follows, we will denote by $u_{m, N, {\rm POD}}(\bm{\mu})$ the HiPOD approximation for the HiMod solution $u_m(\bm{\mu})$
associated with vector ${\mathbf{u}}_{m, N, {\rm POD}}(\bm{\mu})$.
As known, the bottleneck of the projection approach lies in the assembly of $A_{m} (\bm{\mu})$ and $\mathbf{f}_{m} (\bm{\mu})$.
An efficient assembly can be obtained under an affine parameter dependence hypothesis. This requirement will be accomplished in the considered test cases.
Alternative procedures, such as the empirical interpolation method, are adopted in more complex cases~\cite{Rozzabook}.

\subsubsection{Numerical assessment}\label{ADRPOD}
We apply the HiPOD procedure to the test case in Section~\ref{sec:himodADR}.
We recall that the high fidelity solution is provided, in such a case, by a HiMod approximation. 

We identify the parameter $\bm{\mu}$ in \eqref{ParamFullProblem} with the vector ${\bm \mu}=[\nu, b_{x} , b_{y}, \sigma]^T$, which collects some data of problem \eqref{ADRfull}. Concerning the parameter domain, we pick two ranges characterized by a significantly different amplitude, i.e.,
\begin{equation}
D_1 = [1, 100]^4, \quad
D_2 = [1, 10] \times [15, 25] \times [70, 80] \times [20, 30].
\label{eq:adrranges}
\end{equation}
In both the cases, we randomly select $100$ different samples, so that 
$S = \{{\bm \mu}^{(1)}, {\bm \mu}^{(2)}, \ldots, {\bm \mu}^{(100)}\}$. 

\begin{figure}[tb]
\centering
\includegraphics[height=0.37\textwidth]{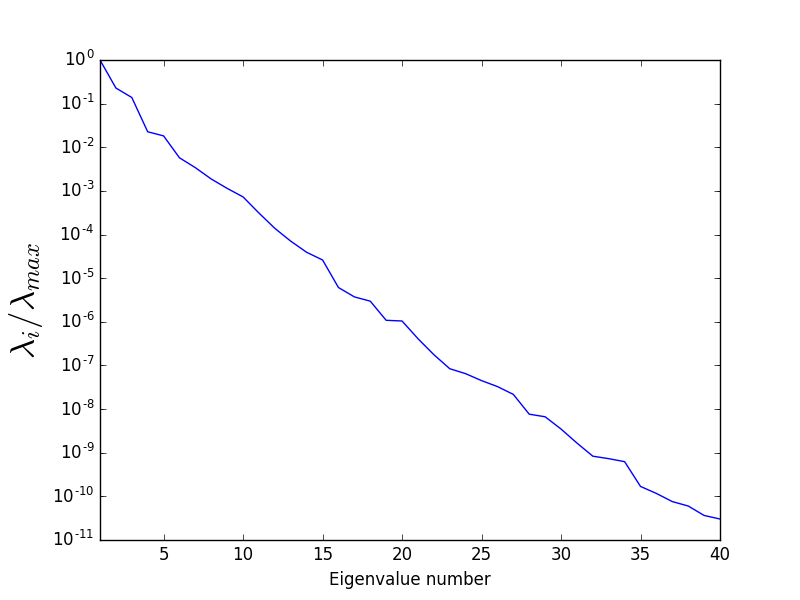} \hfill
\includegraphics[height=0.37\textwidth]{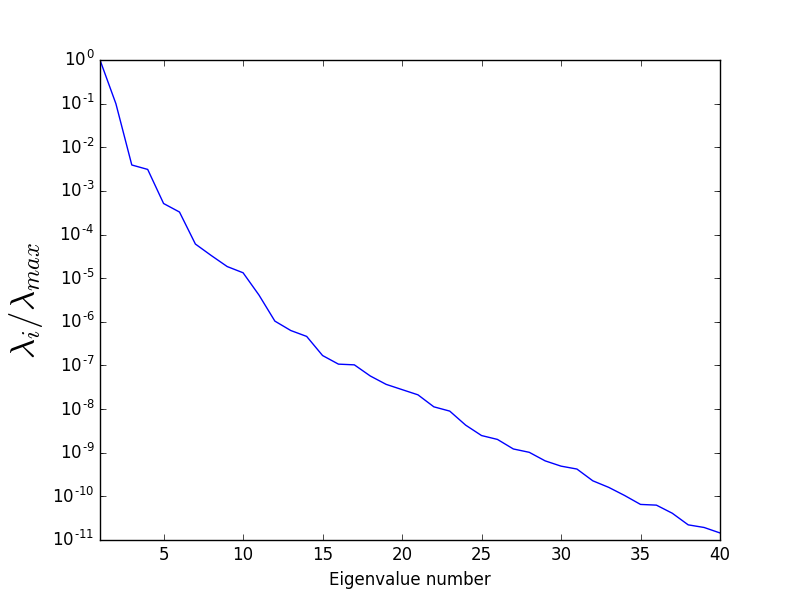}
\caption{ADR test case: eigenvalue trend of the correlation matrix for $D=D_1$ (left) and $D=D_2$ (right), both for $M=100$.}\label{EigErr}
\end{figure}

During the offline phase, we hierarchically reduce the corresponding $100$ ADR problems by employing the same HiMod discretization as in Figure~\ref{adr_f1}, right. \\
The number $N$ of POD basis functions is picked by analyzing the spectrum of the correlation matrix $C$ (see Figure~\ref{EigErr}). The eigenvalues quickly decrease. For the sake of comparison, we select $N=20$ for both ranges. The corresponding eigenvalue, normalized to the maximum one, is $O(10^{-6})$ and $O(10^{-7})$ for $D_1$ and $D_2$, respectively.

Then, we run the online phase to approximate the HiMod solution associated with the parameter
${\bm \mu}=[5, 20, 75, 25]^T$, i.e., the solution in Figure~\ref{adr_f1}, right.

\begin{figure}[tbh]
\centering
\includegraphics[height=0.24\textwidth]{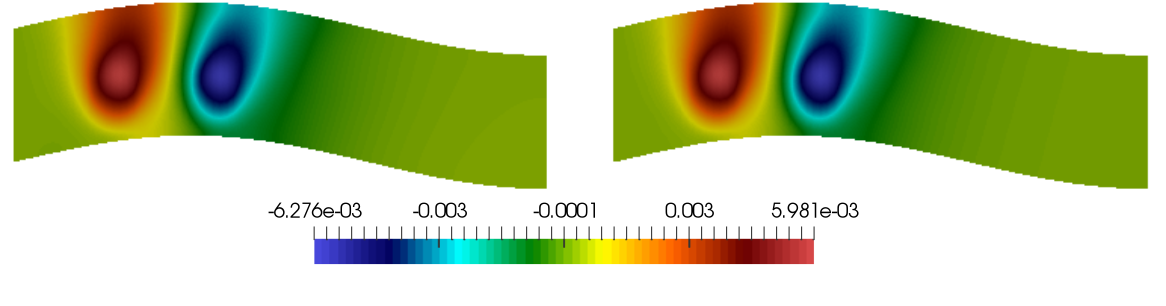}
\caption{ADR test case: HiPOD approximation associated with $D_1$ (left) and $D_2$ (right).}\label{POD_ADR}
\end{figure}

\begin{figure}[tbh]
\centering
\includegraphics[height=0.23\textwidth]{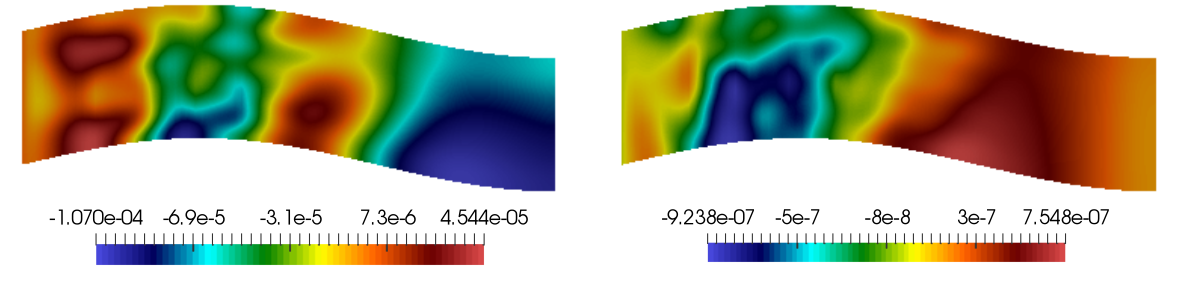}\hfil
\caption{ADR test case: HiPOD modeling error associated with $D_1$ (left) and $D_2$ (right).}\label{PODErrors}
\end{figure}

By comparing the contour plot of the HiPOD approximation in Figure~\ref{POD_ADR} (left for $D_1$, right for $D_2$)
with the HiMod discretization in Figure~\ref{adr_f1}, right, we recognize that 
the global trend of the HiMod solution is correctly detected by both HiPOD solutions. As shown in Figure~\ref{PODErrors},
a more quantitative investigation based on
the distribution of the error $u_{8}(\bm{\mu})-u_{8, 20, {\rm POD}}(\bm{\mu})$ shows that
the solution associated with the smallest parameter domain is, as expected,
more accurate (by about two orders of magnitude) with respect to the approximation obtained when dealing with $D_1$.
The highest accuracy is particularly evident in correspondence with the outflow boundary. 
In Section~\ref{compariamo} we provide a further error analysis for this test case, based on a random sampling.

\subsection{HiPOD reduction for the Stokes equations}
We generalize problem \eqref{FullProblemStokes1} to a parameter dependent setting as follows:
find ${\bm{u}(\bm{\mu})} \in {\bm{V}}, p(\bm{\mu}) \in Q$ such that
\begin{equation}
\begin{cases}
a({\bm{u}}(\bm{\mu}), {\bm{v}}; \bm{\mu}) + b({\bm{v}}, p(\bm{\mu}); \bm{\mu}) = F({\bm{v}}; \bm{\mu}) & \forall \bm{v} \in V,\\
b({\bm{u}}(\bm{\mu}), q; \bm{\mu}) = 0, & \forall q \in Q.
\end{cases}
\label{FullProblemStokes}
\end{equation}

\subsubsection{The offline phase}
As in Section~\ref{OfflinePhase}, we assume a training set $S$ consisting of $M$ different parameters.
According to a \emph{segregated} approach, we generate the POD basis 
for the velocity and for the pressure independently. In~\cite{Baroli} it has been shown 
that a segregated procedure is more effective compared with a monolithic approach, where a unique POD basis (for both velocity and pressure) is built. Moreover, we resort to a vector-valued POD basis for the velocity, in contrast to what has been done in Section~\ref{StokesPODProblem}, where the same (scalar) modal basis is adopted for each component of the velocity. This choice is consistent with standard reduced order modeling techniques in a finite element framework~\cite{ballarin2015supremizer,rozza2007stability}.
Thus, we assemble two distinct response matrices, $U_{\bm{u}}\in \R^{dm_{\bm{u}} N_{h,{\bm{u}}}\times M}$ for the velocity 
and $U_p\in \R^{m_{p} N_{h,p}\times M}$ for the pressure.
Then, by mimicking the scalar case, we compute the correlation matrices associated with $U_{\bm{u}}$
and $U_p$, and we retain the first $N_{\bm{u}}$ and $N_p$ eigenvectors for the velocity and for the pressure, respectively.
This leads us to identify the reduced order spaces ${\bm{V}}_{m_{\bm{u}}, N_{\bm{u}}}$ and 
$Q_{m_p, N_p}$, together with the corresponding matrices $\Upsilon_{m_{\bm{u}}, N_{\bm{u}}}$ and $\Pi_{m_p, N_p}$ collecting, by column, the POD basis functions for velocity and pressure, respectively.

The basic HiPOD procedure is here modified to take into account the stability issue characterizing the approximation provided by a
projection of the Stokes equations. Actually, it turns out that 
even though the solutions involved in the offline phase  
are inf-sup stable, this does not guarantee a priori the inf-sup condition to the reduced space,
with the possible generation of spurious pressure modes.
Following~\cite{ballarin2015supremizer,rozza2007stability}, to recover the inf-sup property for the POD approximation, 
we enrich the velocity space ${\bm{V}}_{m_{\bm{u}}, N_{\bm{u}}}$ with the so-called supremizer solutions. 
In particular, to preserve the offline/online paradigm, we properly modify the procedure proposed in~\cite{ballarin2015supremizer}.
Let ${\mathbf{p}_{{m_p}}}(\bm{\mu}^{(i)})$ be the $i$th column of matrix $U_p$ for $i=1, \ldots, M$. We solve the additional HiMod systems
\begin{equation}
X_{m_{\bm{u}}, \bm{u}}\, {\mathbf{s}_{m_{\bf u}}}(\bm{\mu}^{(i)}) =  B^T_{m_p,m_{\bf u}}(\bm{\mu}^{(i)}){\mathbf{p}_{{m_p}}}(\bm{\mu}^{(i)})
\label{eq:supr}
\end{equation} 
for $=1, \ldots, M$, thus obtaining the $M$ supremizer solutions
${\mathbf{s}_{m_{\bf u}}}(\bm{\mu}^{(i)})\in \R^{dm_{\bm{u}} N_{h,{\bm{u}}}}$.
Here $X_{m_{\bm{u}}, \bm{u}} \in \R^{dm_{\bm{u}} N_{h,{\bm{u}}} \times dm_{\bm{u}} N_{h,{\bm{u}}}}$ denotes the HiMod matrix associated with the inner product in 
${\bm{V}}_{m_{\bm{u}}}$ (so that we use the same modal basis for both velocity and supremizers), while 
$B_{m_p,m_{\bf u}}(\bm{\mu})$ encodes the HiMod discretization of the bilinear form $b({\bm{v}}, p(\bm{\mu}); \bm{\mu})$.
Successively, we assemble the response matrix $U_{\bm{s}}\in \R^{dm_{\bm{u}} N_{h,{\bm{u}}}\times M}$ associated with the supremizers
together with the corresponding correlation matrix,
and we build the matrix $\Xi_{m_{\bm{u}}, N_{\bm{s}}}$ collecting the first $N_{\bm{s}}$ POD
supremizer basis functions, with $N_{\bm{s}} < M$. Finally, we define the matrix
\begin{equation*}
\Phi_{m_{\bm{u}}, N_{\bm{u}} + N_{\bm{s}}} = [\Upsilon_{m_{\bm{u}}, N_{\bm{u}}}, \Xi_{m_{\bm{u}}, N_{\bm{s}}}] \in \R^{dm_{\bm{u}} N_{h,{\bm{u}}}\times (N_{\bm{u}} + N_{\bm{s}})}
\end{equation*}
and the enriched velocity space ${\bm{V}}_{m_{\bm{u}}, N_{\bm{u}} + N_{\bm{s}}}$ spanned by the columns of
$\Phi_{m_{\bm{u}}, N_{\bm{u}} + N_{\bm{s}}}$. 
Throughout the paper, we will refer to ${\bm{V}}_{m_{\bm{u}}, N_{\bm{u}} + N_{\bm{s}}}$
simply as the reduced velocity space. Furthermore, we will always assume $N_{\bm{u}} = N_{\bm{s}} = N_p = N$.

\subsubsection{The online phase}\label{OnlinePhaseStokes}
We extend here the procedure introduced in Section~\ref{OnlinePhase}. For any $\bm{\mu} \in D$, with $\bm{\mu}
\ne \bm{\mu}^{(i)}$ and $i=1, \ldots, M$, rather than solving the corresponding HiMod system
\begin{equation}
\left[\begin{array}{ll} A_{m_{\bm{u}}}(\bm{\mu}) & B^T_{m_p,m_{\bf u}}(\bm{\mu})\\ B_{m_p,m_{\bf u}}(\bm{\mu}) & 0 \end{array}\right]
\left[\begin{array}{c} \mathbf{u}_{m_{\bm{u}}}(\bm{\mu}) \\ \mathbf{p}_{m_p}(\bm{\mu}) \end{array}\right] 
= \left[\begin{array}{c} \mathbf{f}_{m_{\bm{u}}}(\bm{\mu}) \\ \bm{0} \end{array}\right],
\label{HiModProbStokes_param}
\end{equation}
we rely on the reduced system
\begin{equation}
\left[\begin{array}{ll} A_{m_{\bm{u}}, N}(\bm{\mu}) & B^T_{m_p,m_{\bf u}, N}(\bm{\mu})\\ B_{m_p,m_{\bf u}, N}(\bm{\mu}) & 0 \end{array}\right]
\left[\begin{array}{c} \mathbf{u}_{m_{\bm{u}}, 2N}(\bm{\mu}) \\ \mathbf{p}_{m_p, N}(\bm{\mu}) \end{array}\right] 
= \left[\begin{array}{c} \mathbf{f}_{m_{\bm{u}}, N}(\bm{\mu}) \\ \bm{0} \end{array}\right],
\label{pro_POD_Stokes}
\end{equation}
where ${\mathbf{u}}_{m_{\bm{u}}, 2N}(\bm{\mu}) \in \R^{2N}$ and ${\mathbf{p}}_{m_p, N}(\bm{\mu}) 
\in \R^{N}$ denote the POD reduced approximations for the velocity and the pressure, 
$\bm{0}$ is the null vector in $\R^N$, and we assume that matrices
$A_{m_{\bm{u}}, N}(\bm{\mu}) = \Phi_{m_{\bm{u}}, 2N}^T A_{m_{\bm{u}}}(\bm{\mu}) \Phi_{m_{\bm{u}}, 2N}
\in \mathbb \R^{2N \times 2N}$, 
$B_{m_p,m_{\bf u}, N}(\bm{\mu}) = \Pi_{m_{p}, N}^T B_{m_p,m_{\bf u}}(\bm{\mu}) \Phi_{m_{\bm{u}}, 2N}
\in \mathbb \R^{N \times 2N}$ and the vector
$\mathbf{f}_{m_{\bm{u}}, N}(\bm{\mu}) = \Phi_{m_{\bm{u}}, 2N}^T \mathbf{f}_{m_{\bm{u}}}(\bm{\mu})
\in \mathbb \R^{2N}$
can be efficiently assembled owing to affine parameter dependence.
Finally, POD solutions ${\mathbf{u}}_{m_{\bm{u}}, 2N}(\bm{\mu})$ and ${\mathbf{p}}_{m_p, N}(\bm{\mu})$
are projected back to the HiMod space, to yield the approximations
\begin{equation*}
\mathbf{u}_{m_{\bm{u}}}(\bm{\mu}) \approx \Phi_{m_{\bm{u}}, 2N} {\mathbf{u}}_{m_{\bm{u}}, 2N}(\bm{\mu})
:= {\mathbf{u}}_{m_{\bm{u}}, 2N, {\rm POD}}(\bm{\mu}), \quad
\mathbf{p}_{m_p}(\bm{\mu})\approx \Pi_{m_{p}, N} {\mathbf{p}}_{m_p, N}(\bm{\mu})
:= {\mathbf{p}}_{m_p, N, {\rm POD}}(\bm{\mu})
\end{equation*} 
for the HiMod velocity and pressure in \eqref{HiModProbStokes_param}.
The HiPOD approximations for the HiMod solutions $\bm{u}_{m_{\bm{u}}}(\bm{\mu})$ and $p_{m_p}(\bm{\mu})$
will be denoted in what follows by $\bm{u}_{m_{\bm{u}}, 2N, {\rm POD}}$ and $p_{m_p, N, {\rm POD}}$, respectively.

Concerning the stability of the POD reduced problem, 
for any $\bm{\mu} \in D$, one can numerically compare the inf-sup constant associated with the HiMod discretization,
\begin{equation}\label{IS_HiMod}
\beta_{m_{\bm{u}}, m_p}(\bm{\mu}) =  \inf_{\stackrel{q_{{m_p}}\in Q_{{m_p}},}{q_{{m_p}} \neq 0}} 
\ \sup_{\stackrel{{\bm{v}}_{{m_{\bm{u}}}}\in \bm{V}_{{m_{\bm{u}}}},}{{\bm{v}}_{{m_{\bm{u}}}} \neq \bm{0}}} 
\ \frac{b({\bm{v}}_{{m_{\bm{u}}}}, q_{{m_p}}; \bm{\mu})}{\lVert {\bm{v}}_{{m_{\bm{u}}}} \rVert_{\bm{V}} \lVert {q}_{{m_p}} \rVert_{Q}},
\end{equation}
with the corresponding constant resulting from the POD reduction procedure,
\begin{equation}\label{IS_HiPOD}
\beta_{m_{\bm{u}}, m_p, N, {\rm POD}}(\bm{\mu}) =  \inf_{\stackrel{q_{{m_p, N, {\rm POD}}}\in Q_{{m_p, N}},}{q_{{m_p, N, {\rm POD}}} \neq 0}} \ \sup_{\stackrel{{\bm{v}}_{{m_{\bm{u}}, 2N, {\rm POD}}}\in \bm{V}_{{m_{\bm{u}}, 2N}},}{{\bm{v}}_{{m_{\bm{u}}, 2N, {\rm POD}}} \neq \bm{0}}} \ \frac{b({\bm{v}}_{{m_{\bm{u}}, 2N, {\rm POD}}}, q_{{m_p, N, {\rm POD}}}; \bm{\mu})}{\lVert {\bm{v}}_{{m_{\bm{u}}, 2N, {\rm POD}}} \rVert_{\bm{V}} \lVert {q}_{{m_p, N, {\rm POD}}} \rVert_{Q}},
\end{equation}
where, with an abuse of notation, we have adopted the same symbol for the continuous HiPOD spaces as for the corresponding discrete counterparts.
Practical computations for these constants rely on generalized eigenvalue problems (see, e.g.,~\cite{brezzi2013mixed}). In particular,
we resort to the formulas
\begin{equation*}
\beta_{m_{\bm{u}}, m_p}(\bm{\mu}) = \sqrt{\lambda^{(1)}_{m_{\bm{u}}, m_p}}, \quad \beta_{m_{\bm{u}}, m_p, N, {\rm POD}}(\bm{\mu}) = \sqrt{\lambda^{(1)}_{m_{\bm{u}}, m_p, N}},
\end{equation*}
where $\lambda^{(1)}_{m_{\bm{u}}, m_p}$, $\lambda^{(1)}_{m_{\bm{u}}, m_p, N}$, denote the minimum eigenvalue of
the generalized problems
\begin{equation*}
\left[\begin{array}{ll} {X_{m_{\bm{u}}, \bm{u}}} & B^T_{m_p,m_{\bf u}}(\bm{\mu})\\ B_{m_p,m_{\bf u}}(\bm{\mu}) & 0 \end{array}\right]
\left[\begin{array}{c} \mathbf{v}_{m_{\bm{u}}}(\bm{\mu}) \\ \mathbf{q}_{m_p}(\bm{\mu}) \end{array}\right] 
= -\lambda_{m_{\bm{u}}, m_p} \left[\begin{array}{ll} 0 & 0\\ 0 & {X_{m_p, p}} \end{array}\right]
\left[\begin{array}{c} \mathbf{v}_{m_{\bm{u}}}(\bm{\mu}) \\ \mathbf{q}_{m_p}(\bm{\mu}) \end{array}\right],
\end{equation*}
\begin{equation*}
\left[\begin{array}{ll} {X_{m_{\bm{u}}, \bm{u}, 2N}} & B_{m_p,m_{\bf u}, N}(\bm{\mu})^T\\ B_{m_p,m_{\bf u}, N}(\bm{\mu}) & 0 \end{array}\right]
\left[\begin{array}{c} \mathbf{v}_{m_{{\bm{u}, 2N}}}(\bm{\mu}) \\ \mathbf{q}_{m_p, N}(\bm{\mu}) \end{array}\right] 
= -\lambda_{m_{\bm{u}}, m_p, N} \left[\begin{array}{ll} 0 & 0\\ 0 & X_{m_p, p, N} \end{array}\right]
\left[\begin{array}{c} \mathbf{v}_{m_{\bm{u}, 2N}}(\bm{\mu}) \\ \mathbf{q}_{m_p, N}(\bm{\mu}) \end{array}\right],
\end{equation*}
respectively, where $X_{m_{\bm{u}}, \bm{u}}$ is defined as in \eqref{eq:supr}, 
$X_{m_p, p}\in \R^{m_p N_{h,p} \times m_p N_{h,p}}$ denotes the HiMod matrix associated with the inner product in
$Q_{m_p}$, while $X_{m_{\bm{u}}, \bm{u}, 2N}\in \mathbb \R^{2N\times 2N}$, $X_{m_p, p, N}\in \mathbb \R^{N\times N}$ are the corresponding reduced order matrices, given by 
\begin{equation*}
X_{m_{\bm{u}}, \bm{u}, 2N} = \Phi_{m_{\bm{u}}, 2N}^T X_{m_{\bm{u}}, \bm{u}} \Phi_{m_{\bm{u}}, 2N}, \quad 
X_{m_p, p, N} = \Pi_{m_{p}, N}^T X_{m_p, p} \Pi_{m_{p}, N}, 
\end{equation*}
and where, to simplify the notation, we have removed the subscript ${\rm POD}$
to the HiPOD approximations.

\subsubsection{Numerical assessment}\label{podStoksec}
We refer to the test case in Section~\ref{StokesPODProblem}. We identify the parameter with the vector ${\bm \mu}=[\nu, C_{in}, C_{out}, f_{x}, f_{y}]^T$ varying in the domain
\begin{equation*}
D = [1,10] \times [5,15] \times [0,10] \times [1,10] \times [0,10].
\end{equation*}
We consider a sampling set $S$ consisting of $100$ randomly selected values, 
$\{{\bm \mu}^{(1)}, {\bm \mu}^{(2)}, \ldots, {\bm \mu}^{(100)}\}$. 
Then, we hierarchically reduce problem \eqref{FullProblemStokes1} for each parameter in $S$ 
by preserving the same HiMod discretization as the one adopted in Figures~\ref{stokes_f1}-\ref{stokes_f3}, right.

\begin{figure}[tb]
\centering
\includegraphics[height=0.25\textwidth]{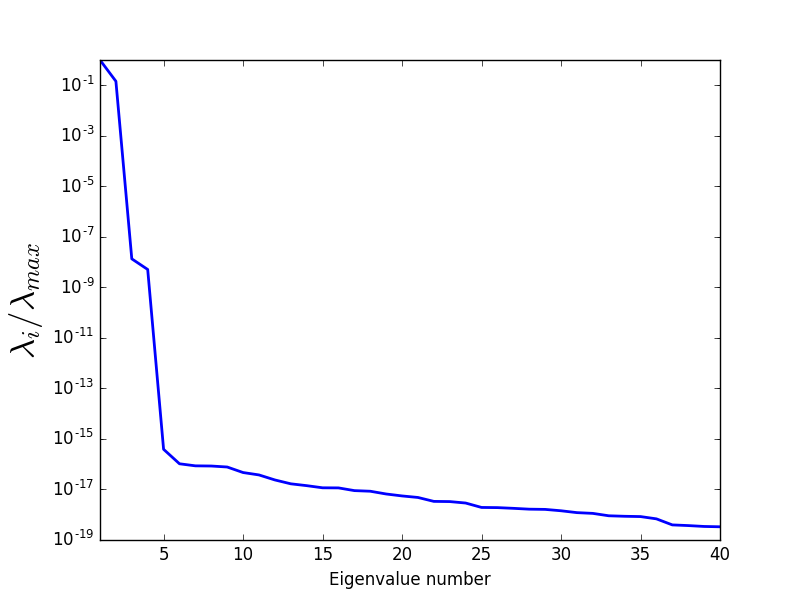}\hfil
\includegraphics[height=0.25\textwidth]{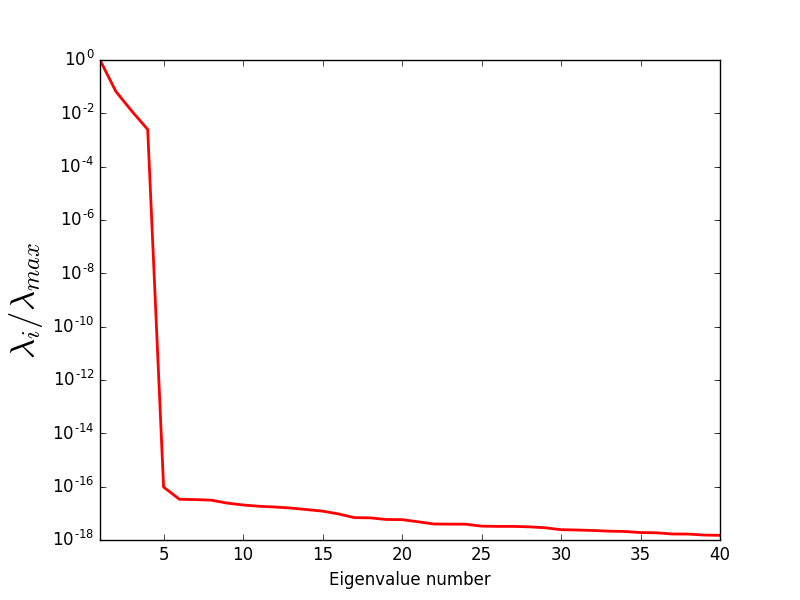}\hfil
\includegraphics[height=0.25\textwidth]{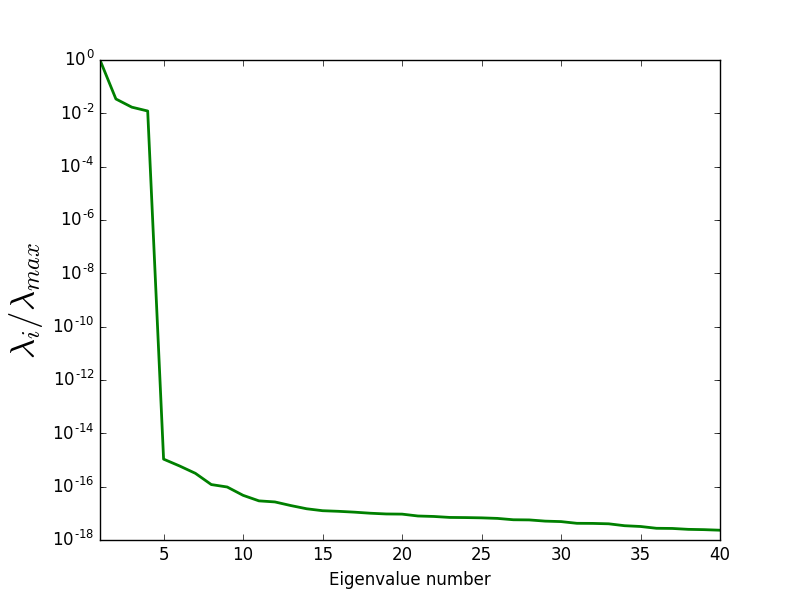}
\caption{Stokes test case: eigenvalue trend of the correlation matrix associated with the velocity (left), the pressure (center), and the supremizers (right).}\label{EigsUPS}
\end{figure}

\begin{figure}[ht!]
\centering
\includegraphics[height=0.21\textwidth]{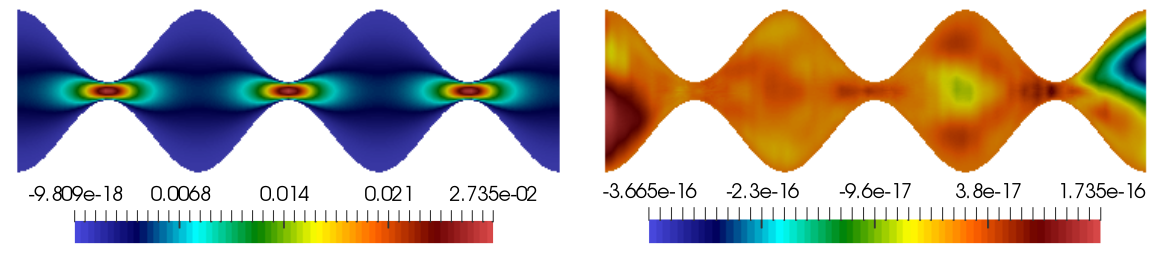}
\caption{Stokes test case: 
HiPOD approximation for the $x$-component of the velocity (left) and associated modeling error (right).}\label{POD_Stokes_x}
\end{figure}

\begin{figure}[ht!]
\centering
\includegraphics[height=0.21\textwidth]{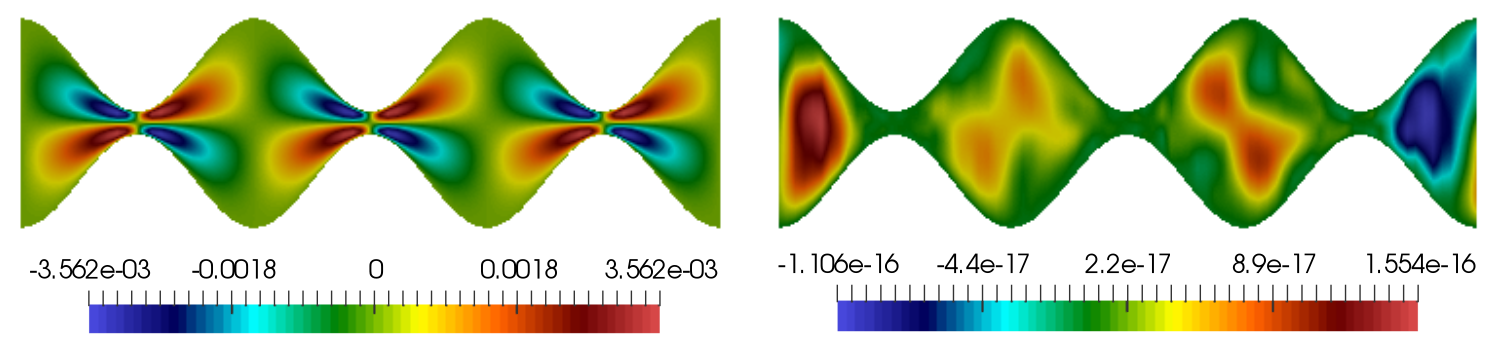}
\caption{Stokes test case: 
HiPOD approximation for the $y$-component of the velocity (left) and associated modeling error (right).}\label{POD_Stokes_y}
\end{figure}

\begin{figure}[ht!]
\centering
\includegraphics[height=0.2\textwidth]{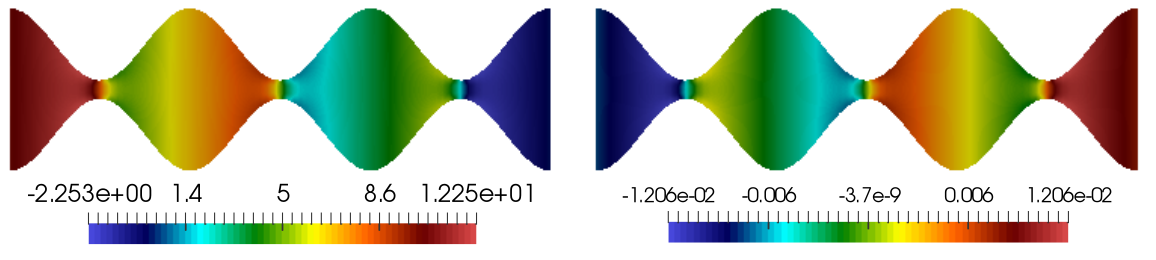}
\caption{Stokes test case: 
HiPOD approximation for pressure (left) and associated modeling error (right).}\label{POD_Stokes_P}
\end{figure}

Figure~\ref{EigsUPS} shows the trend of the eigenvalues of the correlation matrix associated 
with the HiMod velocity, pressure, and supremizers. The drop to the numerical precision occurring
at the fourth eigenvalue in all the plots suggests we set $N=4$.
Indeed, due to the superposition property and since we deal with a linear problem, four independent basis functions are enough to span the space of the solutions to this parametrized Stokes problem.
In the online phase we yield an approximation for the HiMod discretization corresponding to the parameter 
${\bm \mu}=[5, 10, 0, 3, 0]^T$, i.e., for the solution provided in Figures~\ref{stokes_f1}-\ref{stokes_f3}, right.
Figures~\ref{POD_Stokes_x}-\ref{POD_Stokes_P}, left, show the contour plot of the HiPOD approximation for the two components
of the velocity and for the pressure. Comparing the three plots with the corresponding ones 
in Figures~\ref{stokes_f1}-\ref{stokes_f3}, right, we recognize that the selected POD basis 
suffices to provide a reliable approximation, at least
qualitatively. Nevertheless, if we investigate more the details of the distribution of the HiPOD error, we remark that while the two
components of the HiMod velocity are approximated up to machine precision,
this is not the case for the pressure as shown by the contour plots in Figure~\ref{POD_Stokes_P}, right.
This lack of accuracy for the pressure is consistent with what already observed in~\cite{Baroli} and deserves further investigation, for instance, by considering different choices for the supremizers~\cite{rozza2013reduced,abdulle}.

\begin{figure}[ht!]
\centering
\includegraphics[height=0.4\textwidth]{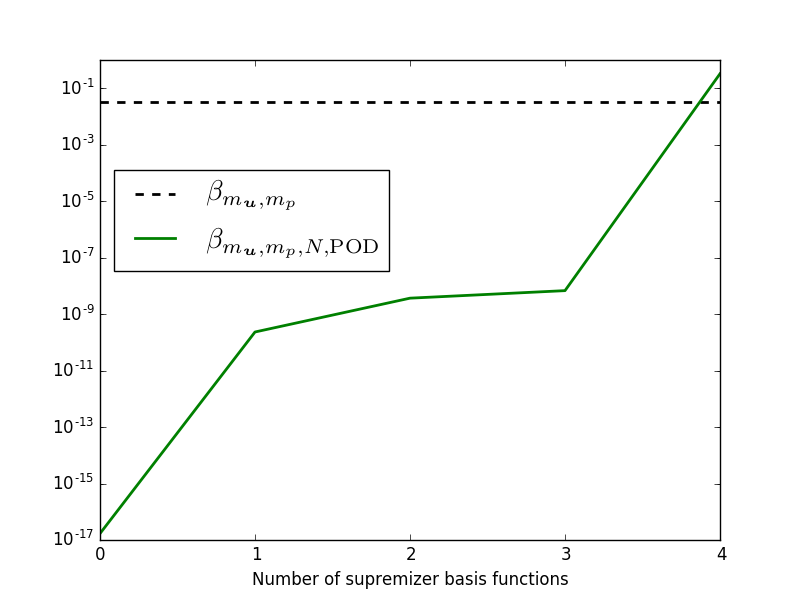}
\caption{Stokes test case: trend of the inf-sup constant $\beta_{m_{\bm{u}}, m_p, N, {\rm POD}}(\bm{\mu})$ as a function of $N_{\bm{s}}$ compared with $\beta_{m_{\bm{u}}, m_p(\bm{\mu})}$.}\label{POD_Stokes_beta}
\end{figure}

Finally, in Figure~\ref{POD_Stokes_beta} we compare the trend of the HiPOD inf-sup constant $\beta_{m_{\bm{u}}, m_p, N, {\rm POD}}(\bm{\mu})$ with $\beta_{m_{\bm{u}}, m_p}(\bm{\mu})$ in \eqref{IS_HiMod}.
For this purpose, we set $N_{\bm{u}} = N_p = 4$ and we make $N_{\bm{s}}$ varying between 0 and 4. To simplify notation, we preserve the same notation as in \eqref{IS_HiPOD} although hypothesis $N_{\bm{u}} = N_p = N_{\bm{s}}$ is here removed.
It is evident that the system is unstable for $N_{\bm{s}} = 0$ so that supremizers are required to recover a reliable pressure, while 
$\beta_{m_{\bm{u}}, m_p, N, {\rm POD}}(\bm{\mu})$ reaches a value comparable with $\beta_{m_{\bm{u}}, m_p}(\bm{\mu})$
when $N_{\bm{s}} = 4$.

\section{The HiRB approach} \label{HiRBSec}
As an alternative to the HiPOD reduction, we introduce a new technique to deal with a
parametrized setting.
We aim at applying an RB approach to HiMod solutions, by relying on a
greedy algorithm during the offline phase.
We refer the interested reader to \cite{rozza2008reduced,Rozzabook,Manzonibook} for an overview of the greedy
algorithm, as well as to~\cite{Maday2002b,binev2011convergence,buffa2012priori,devore2013greedy,hesthaven2014efficient}
for more theoretical insights.
The combination of HiMod with RB justifies the name HiRB adopted to denote the new procedure.

\subsection{HiRB reduction for ADR problems}
Analogously to what was done in Section~\ref{HiPOD_sec}, we first exemplify the HiRB reduction on problem
\eqref{ParamFullProblem}. In particular, we focus on the offline stage since HiPOD and
HiRB essentially resort to the same projection procedure during the online phase.

\subsubsection{The offline phase}\label{OffADRRBSec}
Let $S = \{\bm{\mu}^{(1)}, \bm{\mu}^{(2)}, \ldots, \bm{\mu}^{(M)}\} \subset D^M$ be the training set for the parameter $\bm{\mu}$.
The key idea of the RB method is to generate the reduced space, $\tilde V_{m, N}$ of dimension $N$, by a greedy algorithm, i.e., by adding a single function at a time to the reduced basis~\cite{Rozzabook}.
Let $\tilde V_{m, k}$ denote the reduced space of dimension $k$ known at the $k$th iteration.
Additionally, we assume we have an \emph{error estimator}, $\eta_{m, k}(\bm{\mu})$, for the modeling error associated with
the reduced solution ${u}_{m, k, {\rm RB}}(\bm{\mu})\in \tilde V_{m, k}$, such that
\begin{equation}\label{estim}
\Vert u_{m}(\bm{\mu})-{u}_{m, k, {\rm RB}}(\bm{\mu})  \Vert_{V} \leq \eta_{m, k}(\bm{\mu}) \quad \mbox{with}\ \bm{\mu} \in D,
\end{equation}
with $u_{m}(\bm{\mu})$ the high fidelity HiMod discretization.

Algorithm~\ref{algo:HiRB} itemizes the steps constituting the offline phase of the HiRB method.
First, the
greedy algorithm identifies as a new parameter the value
\begin{equation}\label{stepselection}
\bm{\mu}_g^{(k+1)}  =  \argmax_{\bm{\mu} \in S} \,\eta_{m, k}(\bm{\mu}),
\end{equation}
which corresponds to the most informative HiMod solution not yet included in the reduced space, since it maximizes the discrepancy between
the HiMod space, $V_m$, and the reduced space, $\tilde V_{m, k}$ (step $1$).
The search performed by the greedy algorithm starts from a random choice, $\bm{\mu}_g^{(1)}$, for the parameter and goes on
until a value $\overline{k} \in \mathbb{N}$ is found such that $\max_{\bm{\mu} \in S} \,\eta_{m, \overline{k}}(\bm{\mu}) < \overline{\eta}$,
or the reduced space dimension reaches $N$, with $\overline{\eta}$ a user-defined threshold.
Since the training set $S$ explored in \eqref{stepselection} is finite, a simple enumeration process usually suffices to evaluate the maximum in \eqref{stepselection}, as long as the evaluation of $\eta_{m, k}(\bm{\mu})$ is cheap.
\\
Then, the HiMod approximation,
$\mathbf{u}_{m}(\bm{\mu}_g^{(k+1)}) \in \R^{mN_h}$, is computed (step $2$) and orthonormalized with respect to
the functions already included in the RB basis, collected in the matrix $\tilde \Phi_{m, k} =[\tilde {\bm \varphi}_{1} , \hdots , \tilde {\bm \varphi}_{k}] \in \R^{mN_h\times k}$
(step $3$). This is justified by the fact that solutions at step $2$ can be linearly dependent so that they cannot constitute a basis. We denote the new basis function yielded at step $3$
by $\tilde {\bm{\varphi}}_{k+1}  \in \R^{mN_h}$.
Finally, the RB matrix is extended to include the new basis function, so that we have
\begin{equation*}
\tilde \Phi_{m, k + 1} = [\tilde \Phi_{m, k}, \tilde {\bm{\varphi}}_{k+1}] \in \R^{mN_h\times (k + 1)}
\end{equation*}
(step $4$). This allows us to define the space $\tilde V_{m, k+1}$ of the HiRB approximations associated with $\tilde \Phi_{m, k + 1}$.
\begin{algorithm*}
\normalsize
\caption{HiRB offline phase for ADR problems}\label{algo:HiRB}
\begin{algorithmic}[999]
    \State $1.$ select $\bm{\mu}_g^{(k+1)} = \argmax_{\bm{\mu} \in S} \, \eta_{m,k}(\bm{\mu})$. If $\eta_{m,k}(\bm{\mu}_g^{(k+1)})>\overline{\eta}$, go to $2$, otherwise break;
    \State $2.$ compute $\mathbf{u}_{m}(\bm{\mu}_g^{(k+1)})$;
    \State $3.$ compute the new element, $\tilde {\bm{\varphi}}_{k+1}$, of the basis by orthonormalizing $\mathbf{u}_{m}(\bm{\mu}_g^{(k+1)})$ with respect to $\tilde \Phi_{m, k}$;
    \State $4.$ build the updated RB matrix $\tilde \Phi_{m, k + 1}$ by including $\tilde {\bm{\varphi}}_{k+1}$ as the ($k+1$)th column, and then go back to $1$.
\end{algorithmic}
\end{algorithm*}

Throughout the paper, the HiRB reduced
space eventually yielded by Algorithm~\ref{algo:HiRB} is denoted by
$\tilde V_{m, N}$, with $N$ possibly equal to $\overline{k}$, for $\overline{k}<N$, if the check at step $1$ succeeds at the
$\overline{k}$th iteration.

Finally, the HiRB online phase follows, by mimicking exactly what was performed in Section~\ref{OnlinePhase},
with matrix $ \Phi_{m, N}$ replaced by $\tilde \Phi_{m, N}$. In particular, we denote the HiRB approximation associated with the
vector ${\mathbf{u}}_{m, N, {\rm RB}}(\bm{\mu}) := \tilde \Phi_{m, N} \tilde {\mathbf{u}}_{m, N}(\bm{\mu})$ by $u_{m, N, {\rm RB}}(\bm{\mu})$, with $\tilde {\mathbf{u}}_{m, N}(\bm{\mu})$ the solution of the reduced system corresponding to \eqref{ProjPro}.

The choice of the error estimator, $\eta_{m,k}(\bm{\mu})$, represents a key issue of the RB approach, in particular to ensure the convergence of the greedy algorithm as well as the reliability of the reduced order model.
In general, $\eta_{m,k}(\bm{\mu})$ demands the computation of the reduced solution $\tilde \Phi_{m, k} \tilde {\mathbf{u}}_{m, k}(\bm{\mu})$, so that, at each iteration of Algorithm~\ref{algo:HiRB}, an online phase of dimension $k$ has to be carried out.
For additional details in a standard RB setting, we refer the interested reader, e.g., to~\cite{binev2011convergence,buffa2012priori,devore2013greedy,Rozzabook,rozza2013reduced,rozza2008reduced}. The most common choice for the error estimator relies on the ratio between the dual norm of the weak residual associated with the reduced order solution and a lower bound for the coercivity constant of the high fidelity problem~\cite{Rozzabook}. Thus we have
\begin{equation}\label{etaADR}
\eta_{m,k}(\bm{\mu}) = \frac{\lVert R_{m,k}(v; \bm{\mu}) \rVert_{V'}}{\alpha_{LB}(\bm{\mu})},
\end{equation}
where $V'$ is the dual space of $V$, $R_{m,k}(v; \bm{\mu}) = F(v; \bm{\mu}) - a(u_{m, k, {\rm RB}}(\bm{\mu}), v; \bm{\mu})$ denotes the weak residual of \eqref{ParamFullProblem} associated with the HiRB solution to problem \eqref{ProjPro} for $N := k$, and $\alpha_{LB}(\bm{\mu})$ is a lower bound to the coercivity constant associated with the bilinear form in \eqref{ParamFullProblem}.

The reliability of the chosen estimator is easy to prove. Indeed, we aim to check that
$$\lVert e_{m,k}(\bm{\mu})\rVert_{V} \leq \eta_{m,k}(\bm{\mu})$$,
where $e_{m,k}(\bm{\mu}) = u_m(\bm{\mu}) - {u}_{m, k, {\rm RB}}(\bm{\mu})$.

From \eqref{ParamFullProblem} we have
$$a(e_{m,k}(\bm{\mu}),v;\bm{\mu}) = R_{m,k}(v; \bm{\mu})$$.
Using the definition of dual norm and the coercivity of the bilinear form in \eqref{ParamFullProblem}, we can write
$$a(e_{m,k}(\bm{\mu}),e_{m,k}(\bm{\mu});\bm{\mu}) \leq \lVert R_{m,k}(v; \bm{\mu}) \rVert_{V'} \lVert e_{m,k}(\bm{\mu}) \rVert_V ,$$
$$a(e_{m,k}(\bm{\mu}),e_{m,k}(\bm{\mu});\bm{\mu}) \geq \alpha_{LB} \lVert e_{m,k}(\bm{\mu}) \rVert_V^2 ,$$
i.e.,
$$\lVert e_{m,k}(\bm{\mu}) \rVert_V \leq \frac{\lVert R_{m,k}(v; \bm{\mu}) \rVert_{V'}}{\alpha_{LB}} ,$$
which closes the proof.

In practice, to make the evaluation of the numerator of the error estimator computationally cheap, a Riesz representation property is usually employed, i.e., we look for
$$\hat{r}_{m,k}(\bm{\mu}) \in V : \quad (\hat{r}_{m,k}(\bm{\mu}), v)_V = R_{m,k}(v; \bm{\mu})$$,
where $\hat{r}_{m,k}(\bm{\mu})$ denotes the Riesz representative of $R_{m,k}(\cdot; \bm{\mu})$.
Under the affine parameter dependence assumption of Section \ref{OnlinePhase}, the Riesz representation process considerably simplifies since each term can be represented separately.
For what concerns the denominator of $\eta_{m,k}(\bm{\mu})$, $\alpha_{LB}(\bm{\mu})$ should be cheap to evaluate as well. Nevertheless, standard recipes adopted in the RB literature (such as the successive constraint method \cite{SCM}) need to be suitably modified when adopting HiMod as the high fidelity technique. This represents a topic for a possible future investigation. Here, for simplicity, we set $\alpha_{LB}(\bm{\mu}) := 1$. This choice might not necessarily ensure the reliability of the error estimator (see Section~\ref{compariamo} for a thorough numerical investigation of this issue).

\subsubsection{Numerical assessment}\label{RB_adr_num}
We adopt exactly the same setting as in Section~\ref{ADRPOD}, so that
the parameter ${\bm \mu}$ coincides with the vector $[\nu, b_{x} , b_{y}, \sigma]^T$ and varies
over the ranges, $D_1$ and $D_2$, in \eqref{eq:adrranges}.
Algorithm~\ref{algo:HiRB} is run over a training set $S$ consisting of 100 samples.
Nevertheless, we omit setting the threshold $\overline{\eta}$ at step $1$ and we fix a priori the dimension $N$ of the reduced space to $20$,
also with a view to the comparison performed in Section~\ref{compariamo}.
This leads to hierarchically reducing only $20$ ADR problems, in contrast to $100$ ADR problems with the HiPOD procedure.
Finally, we run the HiRB online phase for ${\bm \mu}=[5, 20, 75, 25]^T$
to approximate the HiMod solution to problem \eqref{ADRfull}.

\begin{figure}[tbh]
\centering
\includegraphics[height=0.24\textwidth]{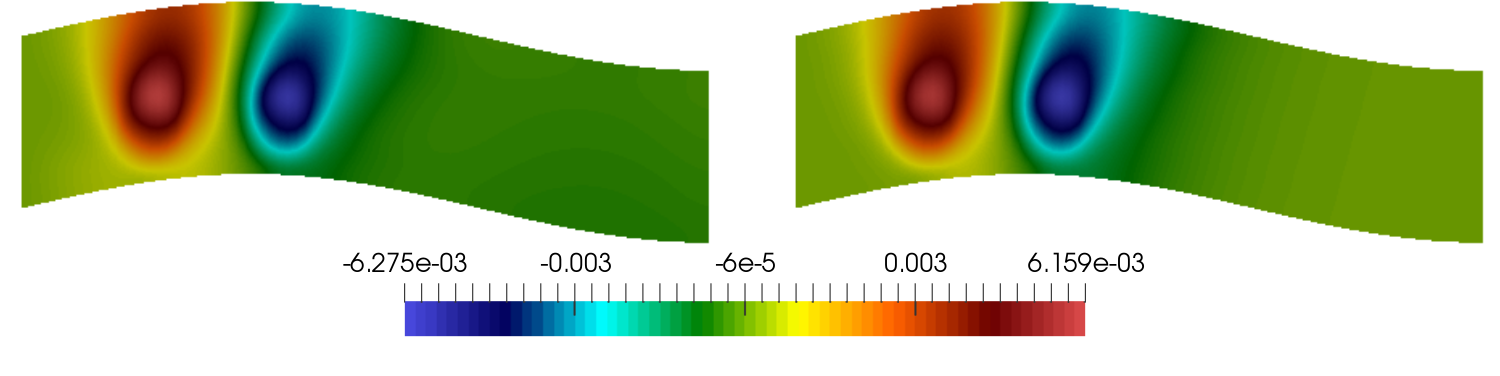}
\caption{ADR test case: HiRB approximation associated with $D_1$ (left) and $D_2$ (right).}\label{RB_ADR}
\end{figure}

\begin{figure}[tbh]
\centering
\includegraphics[height=0.23\textwidth]{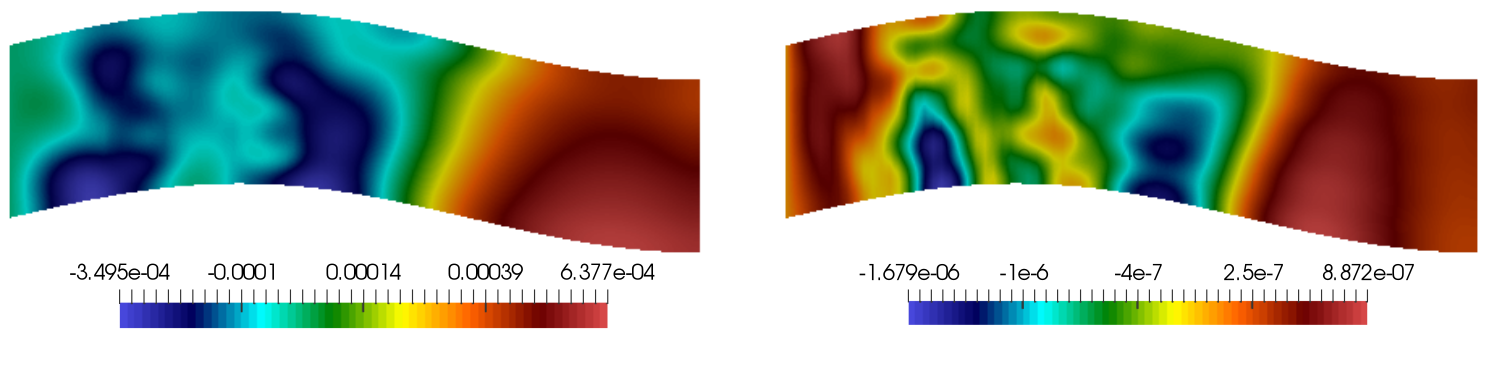}
\caption{ADR test case: HiRB modeling error associated with $D_1$ (left) and $D_2$ (right).}\label{RBErrors}
\end{figure}

Figure~\ref{RB_ADR} shows the contour plot of the HiRB approximation for the two ranges of the parameter ${\bm \mu}$. The qualitative matching between these solutions and the HiMod approximation in Figure~\ref{adr_f1}, right, is good. By analyzing the distribution of the modeling error, $u_{8}(\bm{\mu})-u_{8, 20, {\rm RB}}(\bm{\mu})$,
in Figure~\ref{RBErrors}, it is confirmed that the smaller the parameter range, the higher the accuracy of the HiRB approximation. In particular, the maximum error reduces more than one order when sampling
${\bm \mu}$ in $D_2$.
A cross-comparison with the HiPOD approximations in Figures~\ref{POD_ADR}-\ref{PODErrors} highlights a slightly
higher reliability for the HiRB approach for this test case. A more thorough investigation
in such a direction will be performed in Section~\ref{compariamo}, together with an
error analysis over a random testing.

\subsection{HiRB reduction for the Stokes equations}
We detail the offline step of the HiRB reduction procedure on problem \eqref{FullProblemStokes},
while referring to Section~\ref{OnlinePhaseStokes} for the online phase.

\subsubsection{The offline phase}
In this context, we assume to have an error estimator for both the velocity and the pressure~\cite{Gerner,rozza2013reduced} such that
\begin{equation*}
\Vert {\bm{u}}_{m_{\bm{u}}}(\bm{\mu})- {\bm{u}}_{m_{\bm{u}}, 2k, {\rm RB}}(\bm{\mu})\Vert_{V} +
\Vert p_{m_p}(\bm{\mu})-p_{m_p, k, {\rm RB}}(\bm{\mu})\Vert_{Q}
\leq \eta_{m_{\bm{u}}, 2k, m_p, k}(\bm{\mu}) \quad \mbox{with}\  \bm{\mu} \in D
\end{equation*}
with $\big({\bm{u}}_{m_{\bm{u}}}(\bm{\mu}), p_{m_p}(\bm{\mu})\big)$ the high fidelity HiMod solution pair and
$\big({\bm{u}}_{m_{\bm{u}}, 2k, {\rm RB}}(\bm{\mu}), p_{m_p, k, {\rm RB}}(\bm{\mu})\big)$ the HiRB approximation
belonging to the RB space $\tilde {\bm{V}}_{m_{\bm{u}}, 2k}\times \tilde Q_{m_p, k}$.

As an error estimator, we adopt the quantity $$\eta_{m_{\bm{u}}, 2k, m_p, k}(\bm{\mu}) = \frac{\lVert R_{m_{\bm{u}}, 2k, m_p, k}(v; \bm{\mu}) \rVert_{V'}}{\beta_{LB}(\bm{\mu})} ,$$
where $R_{m_{\bm{u}}, 2k, m_p, k}(v; \bm{\mu})$ is the weak residual of the first equation in \eqref{FullProblemStokes} associated with the HiRB solution when resorting to $m_{\bm{u}}$ and $m_p$ modal basis functions to compute the velocity and the pressure of the HiMod high fidelity space, and $2k$ and $k$ RB functions to evaluate the HiRB velocity and pressure, respectively. Quantity $\beta_{LB}(\bm{\mu})$ is a lower bound for the inf-sup constant associated with the HiMod discretization. Analogously to Section~\ref{OffADRRBSec} the evaluation of this constant is not trivial, especially if we consider that (as discussed in Section \ref{StokesPODProblem}) only empirical criteria are currently available for the choice of compatible HiMod spaces for the velocity and the pressure. For this reason we set $\beta_{LB}(\bm{\mu}) = 1$, while postponing a more rigorous investigation of this issue to a future paper.

Algorithm~\ref{algo:HiRBStokes} details the operations characterizing the HiRB offline phase when applied to the Stokes problem. There are two main differences with respect to Algorithm~\ref{algo:HiRB}, namely (i) we pursue a segregated approach to build the reduced spaces for the velocity and the pressure, (ii) the space for the velocity is enriched via the supremizers.

\begin{algorithm*}
\normalsize
\caption{HiRB offline phase for the Stokes equations}\label{algo:HiRBStokes}
\begin{algorithmic}[999]
    \State $1$ select $\bm{\mu}_g^{(k+1)} = \argmax_{\bm{\mu} \in S} \, \eta_{m_{\bm{u}}, 2k, m_p, k}(\bm{\mu})$. If $\eta_{m_{\bm{u}}, 2k, m_p, k}(\bm{\mu}_g^{(k+1)})>\overline{\eta}$, go to $2$, otherwise break;
    \State $2{\rm a}$ compute the HiMod pair $\big({\mathbf{u}}_{m_{\bm{u}}}(\bm{\mu}_g^{(k+1)}), {\mathbf{p}}_{m_p}(\bm{\mu}_g^{(k+1)})\big)$;
    \State $2{\rm b}$ compute the HiMod supremizer $\mathbf{s}_{m_{\bm{u}}}(\bm{\mu}_g^{(k+1)})$;
    \State $3{\rm a}$ compute the new element, $\tilde {\bm \upsilon}_{k+1}$, of the RB basis for the velocity by orthonormalizing
${\mathbf{u}}_{m_{\bm{u}}}(\bm{\mu}_g^{(k+1)})$ with respect to $\tilde \Upsilon_{m_{\bm{u}}, k}$;
    \State $3{\rm b}$ compute the new element, $\tilde {\bm \pi}_{k+1}$, of the RB basis for the pressure by orthonormalizing ${\mathbf{p}}_{m_p}(\bm{\mu}_g^{(k+1)})$ with respect to $\tilde \Pi_{m_p, k}$;
    \State $3{\rm c}$ compute the new element, $\tilde {\bm \xi}_{k+1}$, of the RB supremizer basis for the velocity by orthonormalizing ${\mathbf{s}}_{m_{\bm{u}}}(\bm{\mu}_g^{(k+1)})$ with respect to $\tilde \Xi_{m_{\bm{u}}, k}$;
    \State $4{\rm a}$ build the updated RB matrix $\tilde \Upsilon_{m_{\bm{u}}, k+1}$ by including $\tilde {\bm \upsilon}_{k+1}$ as the ($k+1$)th column;
        \State $4{\rm b}$ build the updated RB matrix for the pressure $\tilde \Pi_{m_p, k+1}$ by including $\tilde {\bm \pi}_{k+1}$ as the ($k+1$)th column;
            \State $4{\rm c}$ build the updated RB matrix $\tilde \Xi_{m_{\bm{u}}, k+1}$ by including $\tilde {\bm \xi}_{k+1}$ as the ($k+1$)th column;
    \State $4{\rm d}$ build the updated RB matrix for the velocity $\tilde \Phi_{m_{\bm{u}}, 2(k+1)}$, and then go back to $1$.
    \end{algorithmic}
\end{algorithm*}

After the greedy selection on the training set $S$ (step $1$), we solve both the HiMod problem \eqref{HiModProbStokes_param} and the HiMod supremizer equation \eqref{eq:supr} by setting
$\bm{\mu}=\bm{\mu}_g^{(k+1)}$ and $\bm{\mu}^{(i)}=\bm{\mu}_g^{(k+1)}$, respectively,
thus obtaining the HiMod velocity and pressure pair, $\big({\mathbf{u}}_{m_{\bm{u}}}(\bm{\mu}_g^{(k+1)}), {\mathbf{p}}_{m_p}(\bm{\mu}_g^{(k+1)})\big)\in \R^{dm_{\bm{u}} N_{h,{\bm{u}}}}\times
\R^{m_{p} N_{h,p}}$, and the HiMod supremizer, $\mathbf{s}_{m_{\bm{u}}}(\bm{\mu}_g^{(k+1)})\in \R^{dm_{\bm{u}} N_{h,{\bm{u}}}}$ (step $2$).
Then, according to a segregated approach, each HiMod solution is orthonormalized separately, with respect to the corresponding previous basis functions, stored in matrices $\tilde \Upsilon_{m_{\bm{u}}, k}$, $\tilde \Pi_{m_p, k}$, and $\tilde \Xi_{m_{\bm{u}}, k}$, respectively. This yields
the $(k+1)$th RB snapshots, $\tilde {\bm \upsilon}_{k+1}$, $\tilde {\bm \pi}_{k+1}$, $\tilde {\bm \xi}_{k+1}$
(step $3$), which are successively used to enrich the corresponding matrices
(steps $4{\rm a}-4{\rm c}$), so that
\begin{align*}
\tilde \Upsilon_{m_{\bm{u}}, k+1} = [\tilde \Upsilon_{m_{\bm{u}}, k}, \tilde {\bm \upsilon}_{k+1}] \in
\R^{dm_{\bm{u}} N_{h,{\bm{u}}}\times (k+1)}, \\
\tilde \Pi_{m_p, k+1} = [\tilde \Pi_{m_p, k}, \tilde {\bm \pi}_{k+1}] \in \R^{m_{p} N_{h,p}\times (k+1)}, \\
\tilde \Xi_{m_{\bm{u}}, k+1}= [\tilde \Xi_{m_{\bm{u}}, k}, \tilde {\bm \xi}_{k+1}] \in \R^{dm_{\bm{u}} N_{h,{\bm{u}}}\times (k+1)}.
\end{align*}
In particular, matrix $\tilde \Pi_{m_p, k+1}$ allows us to define the $(k+1)$th RB space for the pressure.
The corresponding space for the velocity is the one associated with matrix
$$
\tilde \Phi_{m_{\bm{u}}, 2(k+1)}=[\tilde \Upsilon_{m_{\bm{u}}, k+1}, \tilde \Xi_{m_{\bm{u}}, k+1}]\in
\R^{dm_{\bm{u}} N_{h,{\bm{u}}}\times 2(k+1)},
$$
which is finally built at step $4{\rm d}$.

\subsubsection{Numerical assessment}
We adopt the same numerical setting as in Section~\ref{podStoksec} with the goal of approximating
the HiMod solution in Figures~\ref{stokes_f1}-\ref{stokes_f3}, right, with an RB approach.
Analogously to Section~\ref{RB_adr_num}, we waive
the opportunity to employ the threshold $\overline{\eta}$
in Algorithm~\ref{algo:HiRBStokes}, and we fix the dimension of the reduced spaces to $N=4$
to match the choice in Section~\ref{podStoksec}.

\begin{figure}[ht!]
\centering
\includegraphics[height=0.21\textwidth]{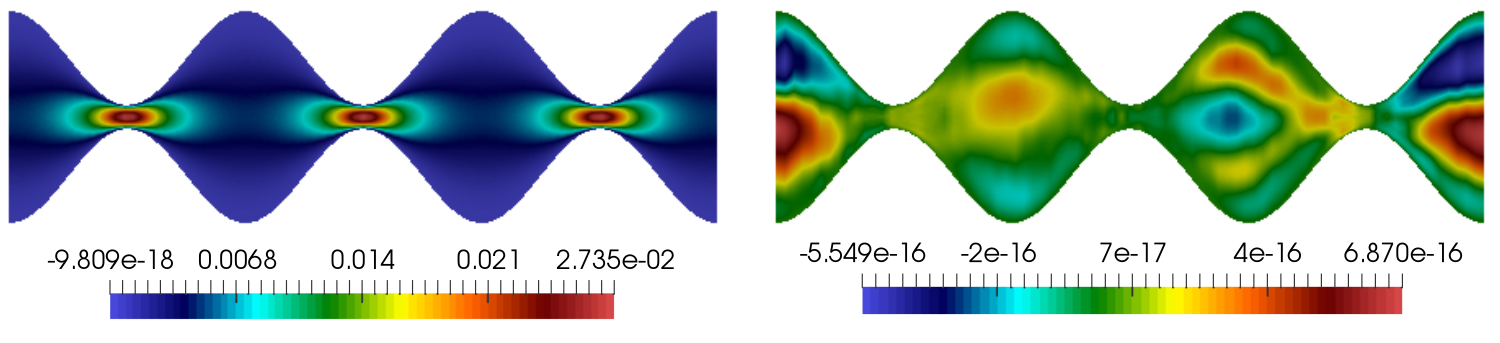}
\caption{Stokes test case:
HiRB approximation for the $x$-component of the velocity (left) and the associated modeling error (right).}\label{RB_Stokes_x}
\end{figure}

\begin{figure}[ht!]
\centering
\includegraphics[height=0.21\textwidth]{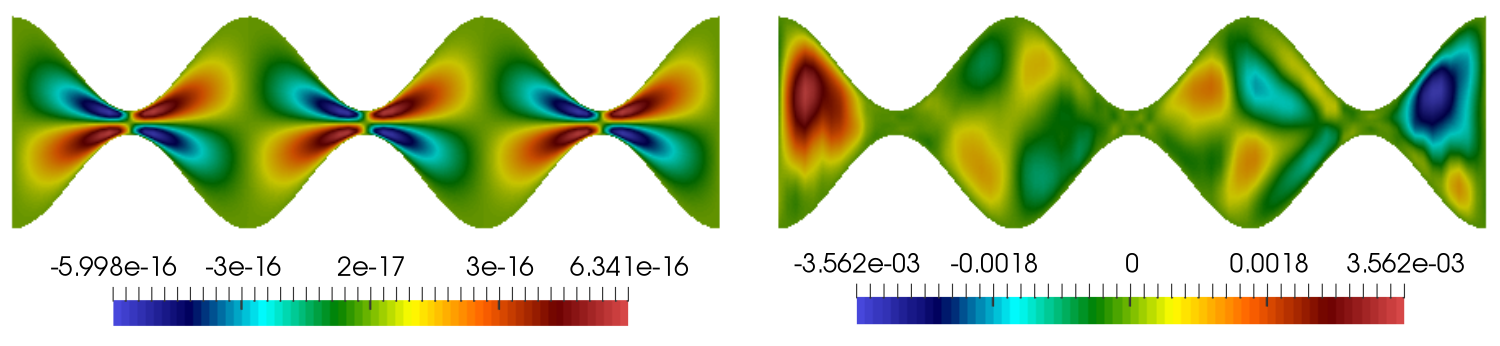}
\caption{Stokes test case:
HiRB approximation for the $y$-component of the velocity (left) and the associated modeling error (right).}\label{RB_Stokes_y}
\end{figure}

\begin{figure}[ht!]
\centering
\includegraphics[height=0.21\textwidth]{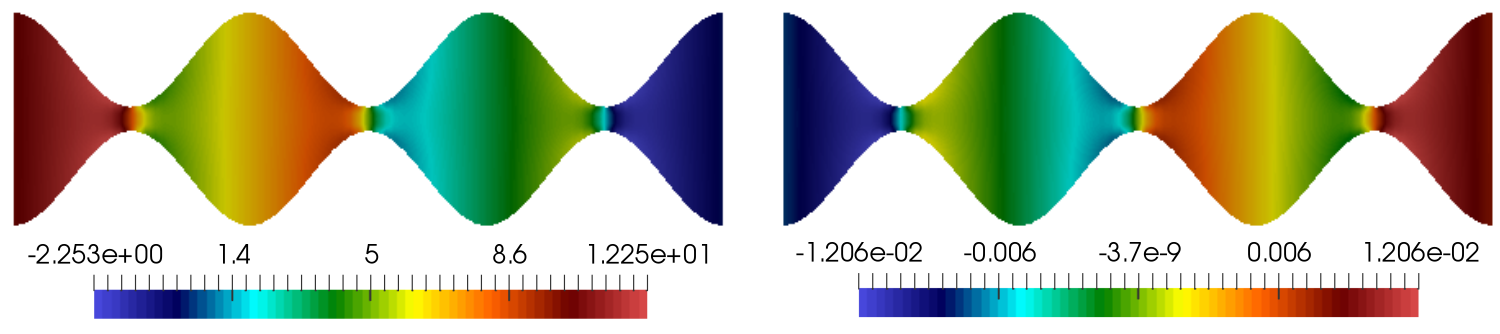}
\caption{Stokes test case:
HiRB approximation for pressure (left) and the associated modeling error (right).}\label{RB_Stokes_P}
\end{figure}

Figures~\ref{RB_Stokes_x}-\ref{RB_Stokes_P}, left, show the contour plot of the HiRB approximation
for the two components of the velocity and for the pressure. The qualitative agreement both with the HiMod
solution and with the HiPOD approximation in Figures~\ref{stokes_f1}-\ref{stokes_f3}, right,
and~\ref{POD_Stokes_x}-\ref{POD_Stokes_P}, left, respectively confirms the reliability of the proposed procedure.
The distribution of the HiRB modeling error in $\Omega$ is provided in
Figures~\ref{RB_Stokes_x}-\ref{RB_Stokes_P}, right. The pressure is the quantity characterized by the worst
accuracy, analogously to what was obtained with the HiPOD approach. Nevertheless, we remark that the HiRB
technique furnishes an approximation of lower quality also for the $y$-component of the velocity
when compared with the HiPOD approximation (notice the difference in terms of order of magnitude for the
corresponding modeling errors in Figure~\ref{RB_Stokes_y}, right, and Figure~\ref{POD_Stokes_y}, right, respectively). Finally, we recognize a more uniform distribution of the modeling error in
Figure~\ref{RB_Stokes_x}, right,
with respect to the corresponding trend of Figure~\ref{POD_Stokes_x}, right. In the former case, the error
is spread across the whole domain, whereas in the latter the error
is mostly confined to the outflow boundary.

\section{HiPOD versus HiRB}\label{compariamo}
This section compares the HiPOD and HiRB techniques. The comparison is carried out 
in terms of three main issues, which are investigated here, separately.

For this purpose, we recall that both HiPOD and HiRB approaches actually carry out a twofold reduction. 
The first one is obtained with the HiMod discretization~\cite{Aletti18,Guzzetti18,Pablo}, while the second reduction
is performed via a projection step during the online phase. The final expectation is the capability to have a reliable approximation for the (full) problem at hand, with a very contained computational effort.

\subsection{Accuracy of the reduced problems}
To compare HiPOD and HiRB in terms of accuracy, 
we plot the average of the associated error over a testing set of $100$ randomly selected parameters 
for both test cases in Sections~\ref{sec:himodADR} and~\ref{StokesPODProblem}.

\begin{figure}[tb]
\centering
\includegraphics[height=0.374\textwidth]{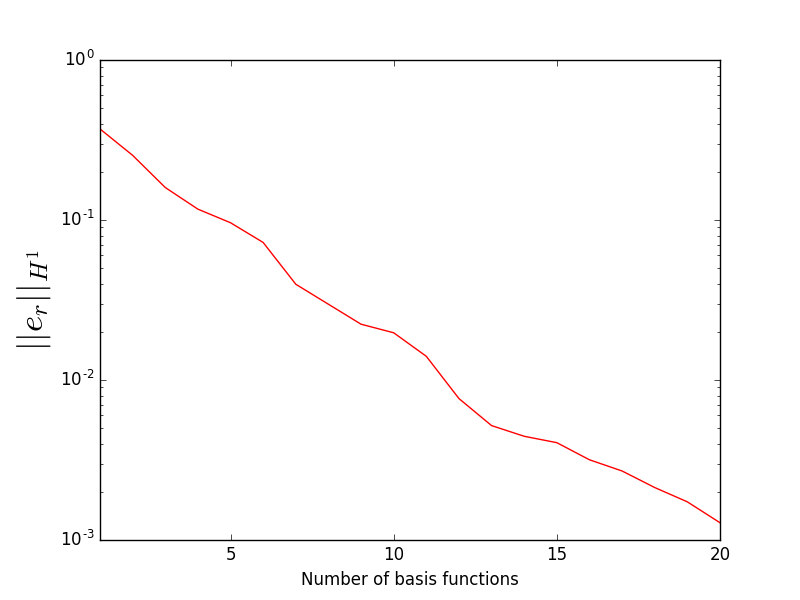}\hfil
\includegraphics[height=0.36\textwidth]{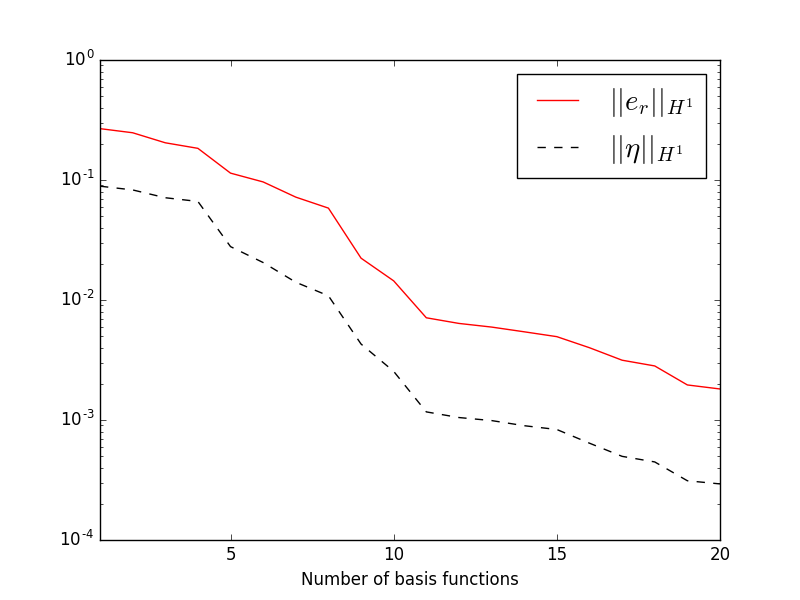}
\caption{ADR test case, parameter range $D_1$: $H^1(\Omega)$-norm of the modeling error associated with the HiPOD (left)
and with the HiRB (right) reduction as a function of $N$.}\label{ErrCompADRLarge}
\end{figure}

\begin{figure}[tb]
\centering
\includegraphics[height=0.374\textwidth]{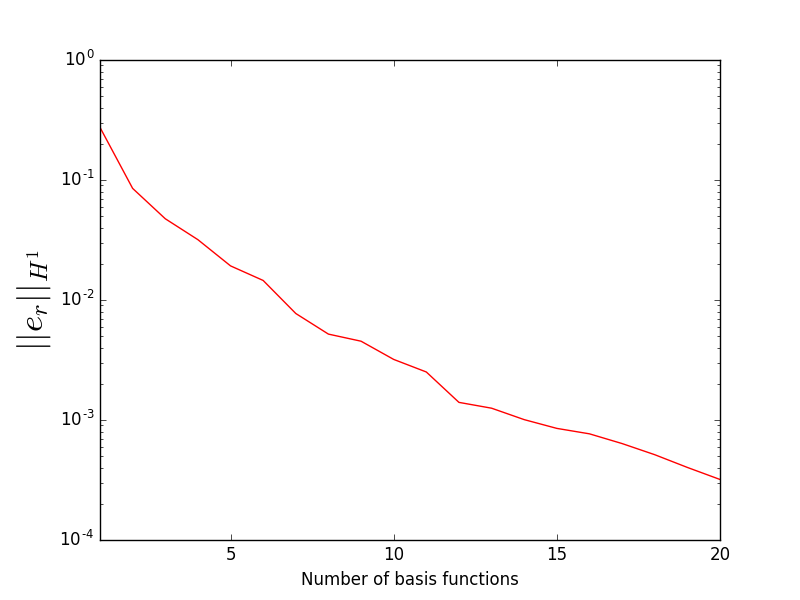}\hfil
\includegraphics[height=0.374\textwidth]{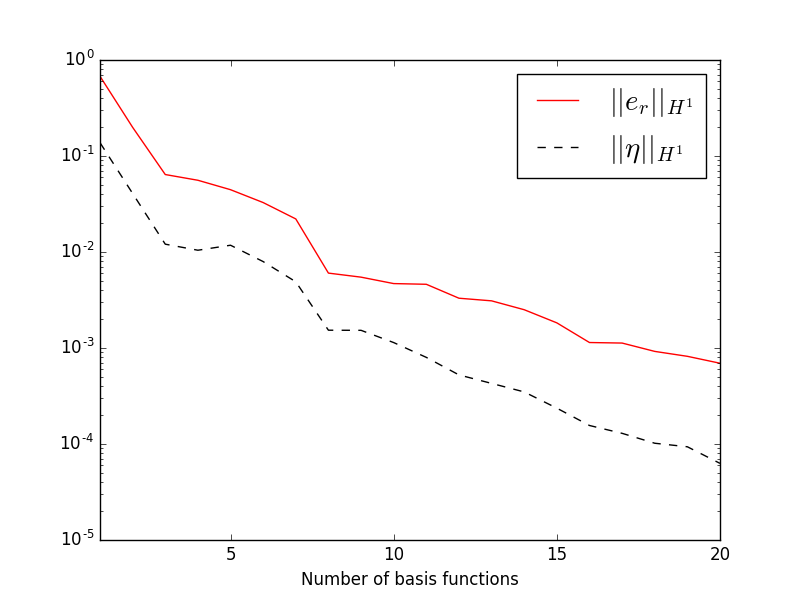}
\caption{ADR test case, parameter range $D_2$: $H^1(\Omega)$-norm of the modeling error associated with the HiPOD (left)
and with the HiRB (right) reduction as a function of $N$.}\label{ErrCompADRRestricted}
\end{figure}

In particular, in Figures~\ref{ErrCompADRLarge} and~\ref{ErrCompADRRestricted} we show the trend of the $H^1(\Omega)$-norm
of the modeling error characterizing the ADR test case and for the two choices of the parameter range in \eqref{eq:adrranges}.
For HiRB approximations, we provide also the trend of the error estimator.
HiPOD exhibits a slightly higher accuracy with respect to HiRB (about half an order of magnitude), for both $D_1$ and $D_2$. 
This is likely related to the adopted estimator which underestimates the exact error of about one order of magnitude
(see Figures~\ref{ErrCompADRLarge} and~\ref{ErrCompADRRestricted}, right).
Actually, $\eta_{m, k}(\bm{\mu})$ is an error indicator rather than an error estimator, since we have set the coercivity constant to one. This might result in a suboptimal greedy selection, although the (monotonic) decreasing trend of the exact error is correctly captured
by $\eta_{m, k}(\bm{\mu})$. Finally, as expected, the reduced order approximation associated with $D_2$  is more accurate for both the procedures.

\begin{figure}[tb]
\centering
\includegraphics[height=0.374\textwidth]{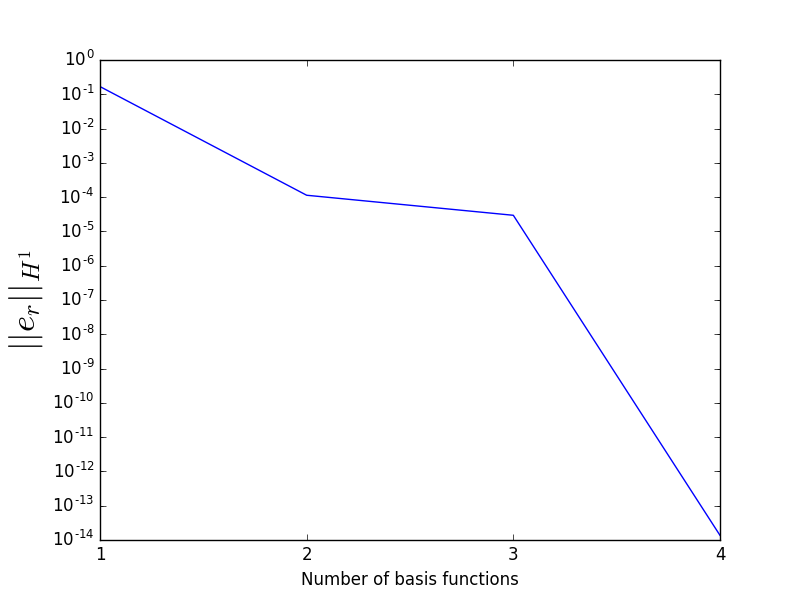}\hfil
\includegraphics[height=0.374\textwidth]{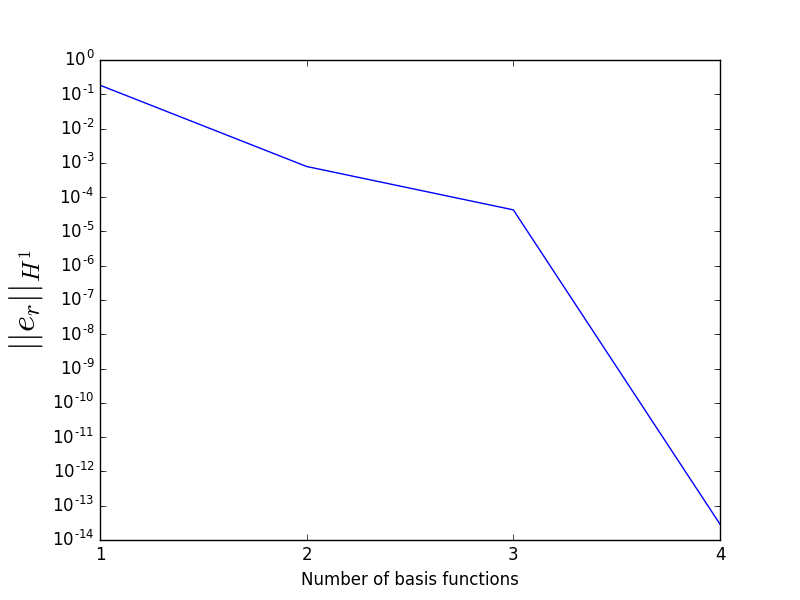}
\caption{Stokes test case: $H^1(\Omega)$-norm of the modeling error associated with the HiPOD (left)
and with the HiRB (right) velocity as a function of $N$.}\label{ErrCompStU}
\end{figure}

\begin{figure}[tb]
\centering
\includegraphics[height=0.374\textwidth]{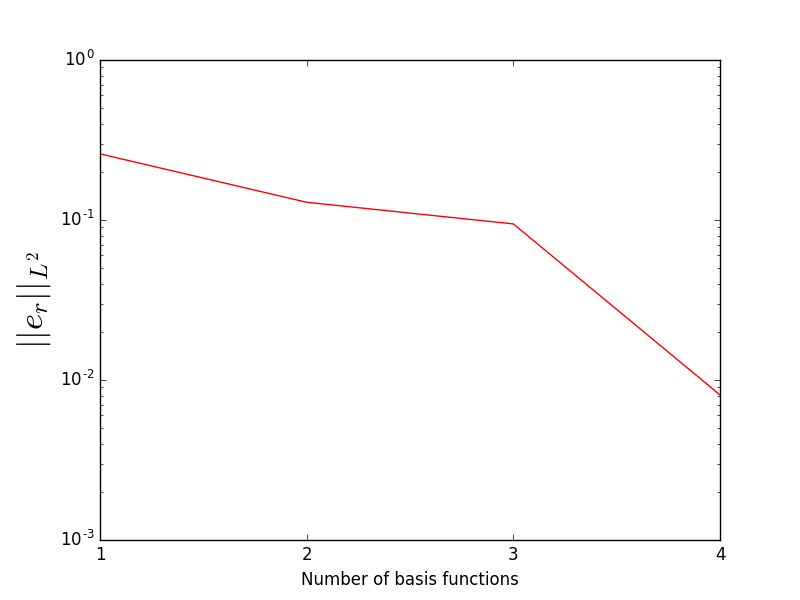}\hfil
\includegraphics[height=0.374\textwidth]{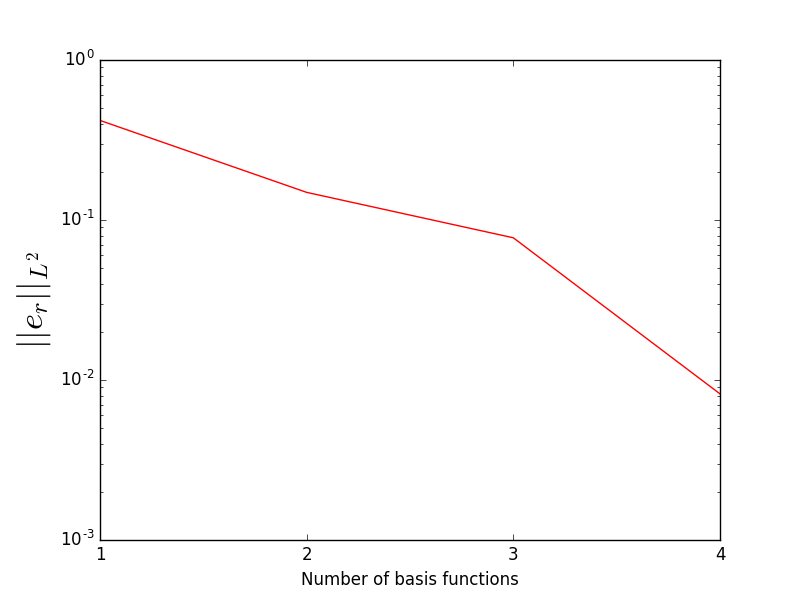}
\caption{Stokes test case: $L^2(\Omega)$-norm of the modeling error associated with the HiPOD (left)
and with the HiRB (right) pressure as a function of $N$.}\label{ErrCompStP}
\end{figure}

The discrepancy between HiPOD and HiRB in terms of accuracy is less evident when considering the Stokes test case (see Figures~\ref{ErrCompStU} and ~\ref{ErrCompStP}).
For $N=4$, the velocity is approximated almost at machine precision, while 
the pressure is characterized by a modeling error of the order of $10^{-2}$ with respect to the $L^2(\Omega)$-norm.
The lower accuracy of the pressure is consistent with what was noted in Figures~\ref{POD_Stokes_P} and~\ref{RB_Stokes_P}.
Possible improvements in such a direction are suggested in Section~\ref{podStoksec}. 
Moreover, the computation of separate error bounds for the velocity and the pressure would be helpful in improving the accuracy, despite requiring further evaluation of stability factors~\cite{Gerner}.

\subsection{Speedup of the reduced problems}
We investigate here the computational effort demanded by the online stage of the HiPOD and HiRB methods. 
We quantify such an effort in terms of CPU time. \footnote{All the simulations are performed on a laptop with an Intel R CoreTM i7 CPU and 4GB RAM.} In particular, we quantify the speedup characterizing the two approaches with the ratio $\tau_{m}(\bm{\mu})/\tau_{m, N}(\bm{\mu})$, where we denote 
the elapsed time associated with the standard HiMod approximation by $\tau_m(\bm{\mu})$ and the time required to solve the corresponding HiPOD or HiRB system by $\tau_{m, N}(\bm{\mu})$. 
A value of speedup greater than one results in a computational gain.

Table~\ref{SpeedAn} gathers the values of this investigation. In order to filter out any dependence on $\bm{\mu}$,
we compute the speedup index over a testing set of $100$ randomly selected parameters.
We observe a large speedup for the Stokes test case and a very mild sensitivity with respect to $N$ for the ADR problem, 
independently of the adopted technique. More in detail, we point out a general lower speedup 
(of about one-third) for the HiRB procedure when compared with HiPOD, in particular for the ADR test case.
This is due to the fact that quantity $\tau_{m, N}(\bm{\mu})$ 
includes also the time elapsed for the evaluation of the error estimator in the HiRB case, whereas this is not required by the HiPOD procedure. Actually, the HiPOD and the HiRB speedups become very similar for the Stokes test case.
\begin{table}
\centering
\begin{tabular}{ccccc}
\toprule
${N}$ & HiPOD-ADR  & HiRB-ADR & HiPOD-Stokes & HiRB-Stokes\\
\midrule
1 & 169.7329 & 68.6832 & 866.3399 & 680.9745\\
2 & 182.3241 & 71.1751 & 977.5722 & 747.6919\\
3 & 178.6220 & 70.5680 & 843.5975 & 693.1631\\
4 & 182.1824 & 70.1519 & 793.0461 & 629.6215\\
\bottomrule
\end{tabular}
\caption{Speedup for HiPOD and HiRB methods applied to the test cases in Sections~\ref{sec:himodADR} and~\ref{StokesPODProblem}. }\label{SpeedAn}
\end{table} 

\subsection{Cost of the offline phase}\label{redb_sec}
We focus now on the offline phase, by comparing the total time required by the HiPOD and HiRB procedures to build 
the reduced basis. 

It is reasonable that, for a fixed dimension, $N$, of the reduced space, 
the HiRB reduction requires less offline time than HiPOD.
Actually, to extract the reduced basis, the HiPOD approach 
computes the HiMod discretization for each of the $M$ parameters in the sampling set, $S$,
and, only a posteriori, compresses such information into a reduced basis of dimension $N$, with $N< M$.
On the contrary, the HiRB method iteratively generates the reduced basis by adding a new basis function at each iteration of the
greedy algorithm. Thus, we compute exactly $N$ HiMod solutions, out of the $M$ possible approximations, associated with the parameters in $S$.

Nevertheless, the offline stage of the HiRB method includes the evaluation of the error estimator during the greedy selection.
A key requirement is that this evaluation is computationally cheap. However,
the construction of the data structures (e.g., higher order tensors~\cite{Rozzabook}) required for this purpose usually entails an additional computational cost, which might dominate the overall offline cost if $M$ is small.
Finally, further less relevant differences between the two methods can be pointed out, such as 
the CPU time required to solve the eigenvalue problem associated with the HiPOD reduction.

\begin{figure}[tb]
\centering
\includegraphics[height=0.5\textwidth]{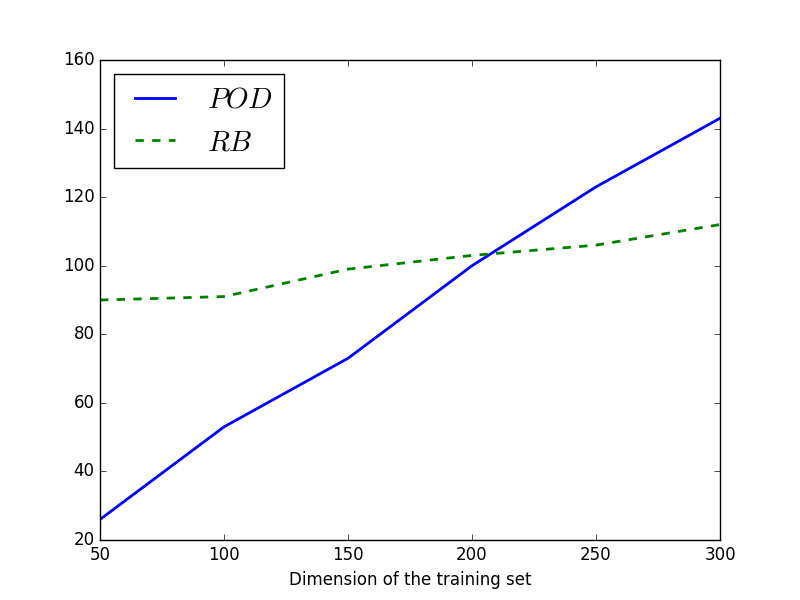}
\caption{ADR test case: comparison between HiPOD and HiRB offline times as a function of $M$.
}\label{TauOff}
\end{figure}

Figure~\ref{TauOff} compares the trend of the total CPU time demanded by the offline stages of the 
HiPOD and HiRB procedures, as a function of the size $M$ of the sampling set, when applied to the test case in Section~\ref{sec:himodADR} and for $N$ set to $20$. 
In agreement with what was noted above, it follows that the 
HiRB training is more expensive than the HiPOD one for small values of $M$.
For instance, for $M=50$, HiPOD takes $25$s, whereas HiRB requires more than $90$s,
most of the time being spent in the setup of the error estimator. Conversely, for large values of $M$,
HiRB demands less time than HiPOD. For example, when $M=300$, HiPOD is more time-consuming 
than HiRB, by requiring $140$s compared with $110$s.

Finally, we remark the different slopes characterizing the two plots in Figure~\ref{TauOff}. The mild slope 
of the HiRB curve confirms that the computational effort required by the evaluation of the error estimator is essentially
independent of the size $M$.
On the contrary, the considerable slope of the HiPOD curve highlights that the computation of the HiMod approximations is not 
negligible and, in general, is heavier than the evaluation of $\eta_{m, k}(\bm{\mu})$.

\section{Conclusions}\label{conclusione}
This work is meant as a first attempt to compare the new reduction techniques HiPOD and HiRB, for
the modeling of parametrized problems. HiPOD has been introduced in~\cite{LupoPasiniPerottoVeneziani,Baroli},
whereas the HiRB approach is proposed here for the first time. The two methods are then compared 
on a benchmark ADR and Stokes problem. Starting from this comparison, we can state that HiRB is better performing than HiPOD
for large training sets, thus turning out to be the ideal tool to tackle, for instance, time demanding fluid dynamics problems. 
As expected, the greedy algorithm allows us to reduce the offline time. The weak point of the HiRB
approach remains the availability of a reliable error estimator. So far,  to simplify the introduction of the new method,
we have adopted an error indicator coinciding with the residual associated with the reduced solution, by completely neglecting the coercivity constant. This rough choice actually 
leads to underestimating the exact error, with a consequent performance loss in terms of speedup and slightly of accuracy with respect to the HiPOD procedure. However,
these conclusions have to be considered as preliminary since we have limited our analysis only to two test cases and, clearly, a more thorough investigation is deserved.

An important issue related to both HiPOD and HiRB has concerned the inf-sup stability which is not necessarily guaranteed for the reduced formulations. To tackle this matter in both cases, supremizer enrichment has been employed with significative improvements.
The pressure approximation for the Stokes equations still demands some amendment for both techniques. 
The proposal of different supremizers likely represents a viable remedy in such a direction.

Among other future developments of possible interest, we cite the generalization to three-dimensional and to nonlinear problems, as 
well as to an unsteady framework. Finally, to certify the reliability of the two methods, a more rigorous investigation of the accuracy characterizing HiPOD and HiRB procedures is desirable, by properly combining HiMod estimates in~\cite{perotto2010hierarchical,Aletti18} with the well-established accuracy results on POD and RB~\cite{Volkwein,Rozzabook}.

\section*{Acknowledments}
This work has been partially supported by the European Union Funding for Research and Innovation, Horizon 2020 Program, in the framework of the European Research Council Executive Agency (H2020 ERC CoG 2015 AROMA-CFD project 681447, ``Advanced Reduced Order Methods with Applications in Computational Fluid Dynamics'', PI Prof. G. Rozza), by the European Union's Horizon 2020 research and innovation program under the Marie Sk\l{}odowska-Curie Actions, grant agreement 872442 (ARIA), and by the research project INdAM 2020, ``Tecniche Numeriche Avanzate per Applicazioni Industriali''. The computations in this work have been performed with the RBniCS library~\cite{rbnics}, developed at SISSA mathLab, which provides an implementation in FEniCS \cite{logg2012automated} of several reduced order modeling techniques. In particular, we acknowledge developers and contributors to both libraries.

\bibliographystyle{plain}
\bibliography{ArticleBibliography}

\end{document}